\documentclass[12pt]{article}
\usepackage[utf8]{inputenc}
\usepackage[all,knot]{xy}
\usepackage{amsmath}
\usepackage{enumerate}
\usepackage{setspace}
\usepackage{graphicx}

\usepackage{tikz}
\usepackage{mathdots}
\usepackage{amsthm}
\usepackage{stmaryrd}
\usepackage{mathtools}
\usepackage{txfonts}

\theoremstyle{plain}
\newtheorem{teor}{Theorem}[subsection]
\newtheorem{prop}[teor]{Proposition}
\newtheorem{cor}[teor]{Corollary}
\newtheorem{lema}[teor]{Lemma}

\theoremstyle{definition}
\newtheorem{defin}[teor]{Definition}

\usepackage{tikz}
\usepackage{dsfont}
\begin{document}
\title{Relative recognition principle}
\author{Renato Vasconcellos Vieira\footnote{renatovv@ime.usp.br\newline
Universidade de São Paulo\newline
The author is greatful for the support provided by CAPES and CNPq through a PhD grant}}
\maketitle

\begin{abstract}
   In this paper the relative recognition principle will be proved. It states that a pair of spaces $(X_o,X_c)$ is weakly equivalent to $(\Omega^N_\text{rel}(\iota:B\hookrightarrow Y),\Omega^N(Y))$ if and only if $(X_o,X_c)$ are grouplike $\overline{\mathcal{SC}^N}$-spaces, where $\overline{\mathcal{SC}^N}$ is any cofibrant resolution of the Swiss-cheese 2-operad $\mathcal{SC}^N$. This principle will be proved for connected $\overline{\mathcal{SC}^N}$-spaces, and also for grouplike $\overline{\mathcal{SC}^N}$-spaces for $2<N\leq\infty$, in the form of an equivalence of homotopy categories. 
\end{abstract}

\section{Introduction}

In this paper the relative recognition principle will be proved. It states that a pair of spaces $(X_o,X_c)$ is weakly equivalent to $(\Omega^N_\text{rel}(\iota:B\hookrightarrow Y),\Omega^N(Y))$ if and only if $(X_o,X_c)$ are grouplike $\overline{\mathcal{SC}^N}$-spaces, where $\overline{\mathcal{SC}^N}$ is any cofibrant resolution of the Swiss-cheese 2-operad $\mathcal{SC}^N$. This principle will be proved for connected $\overline{\mathcal{SC}^N}$-spaces, and also for grouplike $\overline{\mathcal{SC}^N}$-spaces for $2<N\leq\infty$, in the form of an equivalence of homotopy categories. 

In \cite{Ma72} May proved the approximation theorem, which states that there is a natural transformation $\alpha^N$ between the monad $C^N$ associated to the little $N$-cubes operad $\mathcal C^N$ and $\Omega^N\Sigma^N$ that is a weak equivalence on connected spaces. In \cite[2.2]{Ma74} May proved that in general $\alpha^\infty$ is a homological group completion, and in \cite[III.3.3]{CLM76} Cohen extended this result for $1<N<\infty$. The reason it couldn't be extended to the case $N=1$ is that this result requires $C^NX$ and $\Omega^N\Sigma^NX$ to be admissible, and this is only true if $N>1$. In \cite{Ma72} the approximation theorem and the two sided bar construction are used to prove the recognition principle for loop spaces. For an overview including the non-connected cases see \cite{Fr13}. Our proof of the relative recognition principle depends on building a quasifibration $p^N_o:\pi_o\mathcal{SC}^N(X_o,X_c)\rightarrow \mathcal C^{N-1}X_o$ and a map from $p_o$ to a fibration with base space $\Omega^{N-1}\Sigma^{N-1}X_o$ which is $\alpha^{N-1}$ on the base spaces, and then proving that for the conditions of the relative recognition principle these are weak equivalences. From the fact that $\alpha^1$ is not a group completion, it turns out that the proof in this paper does not apply for the non-connected cases when $N=2$, and also for the non-connected cases when $N=1$ for similar reasons. The general case for $N=1$ was proved by Hoefel, Livernet and Stasheff in \cite{HLS16} using the operad of $A_\infty$-actions, which is homotopy equivalent to $\mathcal{SC}^1$.

We will show that for $N<\infty$ there is an equivalence between the homotopy category of $N$-connected relative spaces and the homotopy category of connected $\mathcal{SC}^N$-spaces, and assuming $N>2$ the equivalence in fact holds between the larger homotopy categories of $N-1$-connected relative spaces and of grouplike $\mathcal{SC}^N$-spaces. For $N=\infty$ we will get an equivalence between the homotopy category of 0-connective relative spectra and grouplike $\mathcal{SC}^\infty$-spaces.

The Swiss-cheese 2-operad was introduced by Voronov in \cite{Vo99} as a model of the moduli space of genus-zero Riemann surfaces appearing in the open-closed string theory studied by Zwiebach \cite{Zw98}. Kontsevich used the Swiss-cheese 2-operad in his work on deformation quantization to describe actions of $C^\ast(\mathcal C^N)$-algebras on $C^\ast(\mathcal C^{N-1})$-algebras \cite{Ko99}. Related to Kontsevich's approach to deformation quantization and Zwiebach’s open-closed string field theory Kajiura and Stasheff introduced Open-Closed Homotopy Algebras (OCHA) and Strong Homotopy Leibniz Pairs (SHLP) in \cite{KS06}, which are the algebras over operads that can be obtained from the homology of the Swiss-cheese operad, as has been shown by Hoefel in \cite{HoE07} and by Hoefel and Livernet in \cite{HL12}.

The Swiss-cheese 2-operad itself has been the subject of intense study recently by several authors. Livernet has shown that unlike the little cubes operads the Swiss-cheese 2-operads is not formal \cite{Li15}, and Willwacher has shown that extended Swiss-cheese 2-operads are also not formal \cite{Wi17}. Idrissi has found a model of $\mathcal{SC}^2$ in the category of groupoids \cite{Id17}, and in general Quesney has found combinatorial models for $\mathcal{SC}^N$ in the category of sets and used them to exhibit models for relative loop spaces in dimension 2 \cite{Qu15}.

This paper is organized as follows. In section 2 the structures on the category of topological spaces, sequential spectra, relative spaces and relative spectra that we will need to study relative loop spaces are presented.

In section 3 the theory of monads, operads and operads colored on the ordered set with two elements $2=\{o<c\}$, which we refer to as 2-operads\footnote{Not to be confused with the notion of operads in 2-categories.}, is presented. The little cubes operads and the Swiss-cheese 2-operads and their actions on loop spaces and relative loop spaces are described. We also review some results about the model structure on 2-operads and their algebras.

Readers familiar with operad and model category theory may want to skip sections 2 and 3, with the exception of subsection 2.4 for the definition of the 2-loop spaces and 2-suspension functors and for the unit and counit of their adjunction. All new results are in section 4. First a relative version of a corollary of the approximation theorem is proved. After some technical results on the compatibility of geometric realization with the 2-loop space functors, 2-suspension functors and the monads associated to 2-operads, the relative recognition principle will follow from the above mentioned corollary.

\section{Topological spaces, loop functors and spectra}

We denote by $\texttt{Top}$ the category of compactly generated weakly Housdorff topological spaces. This is a convenient category of topological spaces in the sense that the following hold \cite{StN09}: (1) $\texttt{Top}$ is cartesian closed; (2) $\texttt{Top}$ is bicomplete; (3) Every CW-complex is an object of $\texttt{Top}$; (4) $\texttt{Top}$ is closed under closed subspaces.

The first property allows us to define the $N$-fold loop space functors $\Omega^N$ and guarantees these functors have left adjoints given by the $N$-fold suspension functors $\Sigma^N$. The second allows us to define various constructions as limits and colimits. In particular, the $N$-fold relative loop space functors are defined by applying the $N-1$-fold loop space functors on homotopy fibers of maps, a type of pullback. Note that these (co)limits are not the same as the ones in the category of all topological spaces (see \cite{StN09}), for instance the categorial product in $\texttt{Top}$ has the same underlying set as the one in the category of all topological spaces but a finer topology. The importance of the third and fourth property in the proof of the relative recognition principle presented here is that in the Quillen model structure on topological spaces the cofibrant-fibrant objects (the objects of the homotopy category) are retracts of CW-complexes (retracts of Hausdorff spaces are closed), which means that the homotopy category of $\texttt{Top}$ is isomorphic to the homotopy category of all spaces.
Also, $\texttt{Top}$ has the property that products of CW-complexes are CW-complexes, which isn't true in the category of all topological spaces.

In 2.1 we will present the closed monoidal structure of $\texttt{Top}$ and $\texttt{Top}_\ast$. In 2.2 we present the model structure on $\texttt{Top}$, sketch how homotopy categories are built and present the transfer theorem which allows us to define new model structure on categories through adjunctions. In particular this gives us a model structure on $\texttt{Top}_\ast$. In 2.3 we present the $N$-fold loop space functors and their Quillen left adjoints, the $N$-fold suspension functors, for $1\leq N<\infty$. We also present the category of topological sequential spectra in order to define the $\infty$-fold loop space functor and its Quillen left adjoint, the $\infty$-fold suspension functor. In 2.4 we present the model category of relative spaces and relative sequential spectra in order to define the $N$-fold 2-loop space functors and their left adjoints, the $N$-fold 2-suspension functors. We note that they are only weakly Quillen adjoints.

\subsection{Closed monoidal categories}

A monoidal structure on a category allows us to define a notion of monoid object, which generalizes the notion of monoid from abstract algebra.

\begin{defin}
    A \textit{monoidal category} is a category $\mathcal T$ equipped with a bifunctor $\otimes:\mathcal T\times \mathcal T\rightarrow\mathcal T$ called the \textit{tensor product}, an object $\mathds 1\in\mathcal T$ called the \textit{unit}, a natural isomorphism $\alpha_{X,Y,Z}:(X\otimes Y)\otimes Z\rightarrow X\otimes(Y\otimes Z)$ called the \textit{associator}, a natural isomorphism $\lambda_X:\mathds 1\otimes X\rightarrow X$ called the \textit{left unitor} and a natural isomorphism $\rho_X:X\otimes \mathds 1\rightarrow X$ called the \textit{right unitor},
    satisfying the coherence condition that the diagrams below are commutative:
    $$
    \xygraph{
!{(-1.5,0) }*+{(W\otimes (X\otimes Y))\otimes Z}="a1"
!{(1.5,0) }*+{W\otimes ((X\otimes Y)\otimes Z)}="a2"
!{(2.25,1) }*+{W\otimes(X\otimes (Y\otimes Z))}="a3"
!{(0,1.75) }*+{(W\otimes X)\otimes (Y\otimes Z)}="a4"
!{(-2.25,1) }*+{((W\otimes X)\otimes Y)\otimes Z}="a5"
"a1":@{->}_{\alpha}"a2"
"a2":@{->}_{1\otimes \alpha}"a3"
"a5":@{->}^{\alpha}"a4"
"a4":@{->}^{\alpha}"a3"
"a5":@{->}_{\alpha\otimes 1}"a1"
   }\ \raisebox{1.8cm}{\xymatrix@C=0.3cm{
   (X\otimes\mathds 1)\otimes Y\ar[rr]^{\alpha}\ar[dr]_{\rho\otimes 1}&&X\otimes(\mathds 1\otimes Y)\ar[dl]^{1\otimes\lambda}\\
   &X\otimes Y&
   }}$$
   
   A monoidal category is \textit{cartesian} if the tensor product is given by the categorial product $\times$. In this case the unit is the terminal object.
\end{defin}

\begin{defin}
    A \textit{monoid} in a monoidal category $(\mathcal T,\otimes,\mathds 1)$ is an object $M\in\mathcal T$ equipped with a pair of morphisms $\mu:M\otimes M\rightarrow M$, called the \textit{multiplication}, and $\eta:\mathds 1\rightarrow M$, called the \textit{unit}, such that the diagrams below commute:
    
    $$
        \xymatrix{
            M\otimes M\otimes M\ar[r]^{1\otimes \mu}\ar[d]_{\mu\otimes {1}}&M\otimes M\ar[d]^\mu&&&M\ar[r]^{1\otimes \eta}\ar[dr]_{1}&M\otimes M\ar[d]|\mu&M\ar[l]_{\eta\otimes {1}}\ar[dl]^{1}\\
            M\otimes M\ar[r]_\mu&M&&&&M&
        }
    $$
    
    A \textit{morphism of monoids} is a morphism between the objects that commute with the structural morphisms.
    
    A \textit{comonoid} is dually defined by inverting all arrows in the definition of a monoid.
\end{defin}

\begin{defin}
    A monoidal category is \textit{symmetric} if in addition it is equipped with a natural isomorphism $\tau_{X,Y}:X\otimes Y\rightarrow Y\otimes X$ such that $\tau_{Y\otimes X}\tau_{X\otimes Y}=1_{X\otimes Y}$ and satisfying the coherence condition that the diagrams of the type below are commutative
    $$
        \xymatrix{ 
	& X \otimes (Y \otimes Z) \ar^{\tau}[r] & (Y \otimes Z) \otimes X \ar^{\alpha}[dr] & \\
	(X \otimes Y) \otimes Z \ar^{\alpha}[ur] \ar_{\tau \otimes 1}[dr] & & & Y \otimes (Z \otimes X) \\
	& (Y \otimes X) \otimes Z \ar_{\alpha}[r] & Y \otimes (X \otimes Z) \ar_{1 \otimes \tau}[ur] & 
	}$$
	
	A cartesian monoidal category is always symmetric.
\end{defin}
    
    \begin{defin}
        An \textit{adjunction} is a pair of functors $S:\mathcal T\rightarrow \mathcal A$, called the \textit{left adjoint functor}, and 
        $\Lambda:\mathcal A\rightarrow \mathcal T$, called the \textit{right adjoint functor}, equipped with a pair of natural transformations  $\eta:Id_{\mathcal T}\rightarrow \Lambda S$, called the \textit{unit}, and $\epsilon:S\Lambda\rightarrow Id_{\mathcal A}$, called the \textit{counit}, such that the unit-counit diagrams below are commutative
    $$
    \xymatrix{
        S\ar[r]^{S\eta}\ar@/_0.4cm/[rr]_{1_S}&S\Lambda S\ar[r]^{\epsilon_S}&S
    }\ \ \ \ 
    \xymatrix{
        \Lambda\ar[r]^{\eta_{\Lambda}}\ar@/_0.4cm/[rr]_{1_\Lambda}&\Lambda S\Lambda\ar[r]^{\Lambda\epsilon}&\Lambda
    }
    $$

    We denote adjunctions as $(S\dashv \Lambda)$.
    \end{defin}
    
    \begin{defin}
        A \textit{closed monoidal category} is a monoidal category $\mathcal T$ equipped with a bifunctor $-^-:\mathcal T\times \mathcal T^\text{op}\rightarrow \mathcal T$, called the \textit{internal morphisms functor} or \textit{exponential}, such that for every $X\in\mathcal T$ there is an adjunction $(-\otimes X\dashv-^X)$.
    \end{defin}

The category $\texttt{Top}$ of compactly generated weakly Housdorff topological spaces has a closed cartesian category structure, with the cartesian product $\times$ as tensor product and the mapping space with the compact-open topology as the internal morphisms functor.\cite{StN09}

The category $\texttt{Top}_\ast$ of pointed topological spaces $(X,\ast_X)$ also has a closed symmetric monoidal structure. The monoidal product is given by the smash product
$$
    X\wedge Y:=\frac{X\times Y}{\{\ast_X\}\times Y\cup X\times\{\ast_Y\}}
$$ and the unit is $\mathds S^0$. The internal morphisms functor is the pointed mapping space. Sometimes we denote a pointed space $(X,\ast_X)$ simply as $X$.

Note that the forgetful functor
\begin{align*}
    U:\texttt{Top}_\ast&\rightarrow \texttt{Top}\\
    (X,\ast_X)&\mapsto X 
\end{align*}
is left adjoint to the disjoint base point functor
\begin{align*}
    -_\ast:\texttt{Top}&\rightarrow \texttt{Top}_\ast\\
    X&\mapsto \left(X\coprod \ast,\ast\right) 
\end{align*}

\subsection{Model categories}

Given a locally small category $\mathcal T$ and a class of morphisms $W$ we sometimes want to treat them as if they were isomorphisms. For instance in $\texttt{Top}$ we often want to treat weak homotopy equivalences as if they were isomorphisms, despite the fact they don't necessarily have inverses, not even up to homotopy. If there are classes of morphisms satisfying certain lifting properties compatible with $W$ we can define a new locally small category $\mathcal Ho\mathcal T$, called the homotopy category of $\mathcal T$ with respect to $W$, with the universal property that there is a functor $q:\mathcal T\rightarrow\mathcal Ho\mathcal T$ such that:
\begin{enumerate}
    \item for any functor $F:\mathcal T\rightarrow \mathcal U$ with $Fw$ being an isomorphism if $w\in W$, there is a functor $F_q:\mathcal Ho\mathcal T\rightarrow U$ and a natural isomorphism $F\rightarrow F_q\circ q$;
    
    \item the functor $q^\ast:Fun(\mathcal Ho\mathcal T,\mathcal U)\rightarrow Fun(\mathcal T,\mathcal U)$ is full and faithful.
\end{enumerate}
 
The axioms of model categories gives us general conditions under which the above construction is possible. For a detailed account of model categories see \cite{Hi09,HoM07}, and for the definition used here see \cite{Ri09}. First, we require $W$ to satisfy a few conditions.

\begin{defin}
    A \textit{category with weak equivalences} is a category $\mathcal T$ equipped with a subcategory $W\subset \mathcal T$ such that:
    \begin{enumerate}[1.]
        \item $W$ contains all isomorphisms of $\mathcal T$,
        
        \item $W$ satisfy the \textit{two-out-of-three property}: for $f,g$ any pair of composable morphisms of $\mathcal T$, if two of $\{f,g,fg\}$ are in $W$, then so is the third.
    \end{enumerate}
\end{defin}

We could formally add inverses to the elements of $W$ but this may result in a non-locally small category. Homotopy theory gets around this by identifying two classes of morphisms in $\mathcal T$, called cofibrations and fibrations, that have special lifting properties compatible with $W$.

Given two morphisms $f:X\rightarrow X'$ and $g:Y\rightarrow Y'$ in $\mathcal T$ we say that $f$ has the \textit{left lifting property with regards to $g$}, and equivalently that $g$ has the \textit{right lifting property with regards to $f$}, if for every commutative square
$$
    \xymatrix{
        X\ar[d]_f\ar[r]^{H_0}&Y\ar[d]^g\\
        X'\ar[r]_{H}\ar@{-->}[ur]|{\tilde H}&Y'
    }
$$
a lift $H$ exists that makes the two triangles commute. For a class of morphisms $M$ in $\mathcal T$ we denote the class of all morphisms with the left lifting property with regards to all morphisms in $M$ as ${}^\boxslash M$, and the class of all morphisms with the right lifting property with regards to all morphisms in $M$ as $M^\boxslash$.

In the definition bellow $\mathcal{T}^{\shortrightarrow}$ is the category of morphisms in $\mathcal T$ and $\mathcal{T}^{\shortrightarrow\shortrightarrow}$ the category of composable morphisms in $\mathcal T$.

\begin{defin}
    A \textit{weak factorization system} on a category $\mathcal T$ is a pair $(L,R)$ of classes of morphisms of $\mathcal T$ such that $L={}^\boxslash R$ and $R=L^\boxslash$, equipped with a functor
        \begin{align*}
            \text{Fact}_{L,R}:\mathcal{T}^{\shortrightarrow} &\rightarrow \mathcal{T}^{\shortrightarrow\shortrightarrow}\\
            (f:X\rightarrow Y)&\mapsto (X\xrightarrow{f_L} \text{Fact}_{L,R}(f)\xrightarrow{f_R}Y)
        \end{align*}
        such that $f=f_Rf_L$, $f_L\in L$ and $f_R\in R$.
\end{defin}

\begin{defin}
    A \textit{model structure} on a category $\mathcal T$ is a choice of three distinguished classes of morphisms $(W,C,F)$ in $\mathcal T$, the \textit{weak equivalences} $W$, \textit{cofibrations} $C$ and \textit{fibrations} $F$ satisfying the following conditions
    \begin{enumerate}[1.]
        \item $(\mathcal T,W)$ is a category with weak equivalences;
        
        \item There are functors $\text{Fact}_{C,W\cap F}$ and $\text{Fact}_{C\cap W, F}$ such that $(C,W\cap F)$ and $(C\cap W, F)$ are two weak factorization systems on $\mathcal T$.
    \end{enumerate}
\end{defin}

We refer to morphisms in $C\cap W$ as \textit{trivial cofibrations} and morphisms in $W\cap F$ as \textit{trivial fibrations}.

\begin{defin}
    A \textit{model category} is a bicomplete category $\mathcal T$ equipped with a model structure $(W,C,F)$.
\end{defin}

If $\mathcal T$ is bicomplete it in particular has an initial and final object. Since our primary example of model category will be $\texttt{Top}$ we denote initial objects as $\emptyset$ and final objects as $\ast$. We say that an object $X$ is \textit{cofibrant} if $\emptyset \rightarrow X\in C$, and that it is \textit{fibrant} if $X\rightarrow \ast\in F$.

For $X$ a cofibrant object and $Y$ a fibrant object in $\mathcal T$ there is a well defined equivalence relation on the set of morphisms $\mathcal T(X,Y)$, where morphisms on the same equivalence classes are said to be homotopic (see \cite[7.3, 7.4]{Hi09} for details).

The weak factorization structures of a model category gives us natural ways to associate a weakly equivalent (co)fibrant object to any object in a model category.

\begin{defin}
In a model category $\mathcal T$ the \textit{cofibrant replacement functor} is
\begin{align*}
    \text{Cof}:\mathcal T&\rightarrow\mathcal T\\
    X&\mapsto\text{Fact}_{C,W\cap F}(\emptyset\rightarrow X)
\end{align*}
and there is a natural trivial fibration $\text{cof}:\text{Cof}\rightarrow Id_\mathcal T$ with $\text{cof}_X:=(\emptyset\rightarrow X)_{W\cap F}$.

The \textit{fibrant replacement functor} is
\begin{align*}
    \text{Fib}:\mathcal T&\rightarrow\mathcal T\\
    X&\mapsto\text{Fact}_{C\cap W,F}(X\rightarrow\ast)
\end{align*}
and there is a natural trivial fibration $\text{fib}:Id_\mathcal T\rightarrow \text{Fib}$ with $\text{fib}_X:= (X\rightarrow\ast)_{C\cap W}$.
\end{defin}

In general, we refer to a cofibrant object $\overline X$ equipped with a weak equivalence $\overline X\rightarrow X$ as a \textit{cofibrant resolution} of $X$.

The Whitehead theorem for model categories then states that a morphism between cofibrant-fibrant objects is a weak equivalence if and only if it is a homotopy equivalence. This justifies the following definition.

\begin{defin}
    The \textit{homotopy category} $Ho\mathcal T$ of a model category $\mathcal T$ is the category whose objects are cofibrant-fibrant objects of $\mathcal T$, and morphisms are the homotopy equivalence classes of morphisms between them.
\end{defin}

The functor $q:\mathcal T\rightarrow \mathcal Ho\mathcal T$ defined on objects as $q(X):=\text{Fib}\text{Cof}(X)$ and on morphisms as $q(f):=[\text{Fib}\text{Cof}(f)]$ satisfies the universal property from the beginning of this section.

The Quillen model structure $(W_Q,C_Q,F_Q)$ on $\texttt{Top}$ is given by the following classes of maps \cite{Hi15}:
\begin{itemize}
    \item $f:X\rightarrow Y\in W_Q$ if it is a weak homotopy equivalence, i.e.  $f_\ast:\pi_q(X,x)\rightarrow \pi_q(Y,y)$ is an isomorphism for all $x\in X,\ y\in Y$ and $q\in\mathds N$;
    
    \item $f:X\rightarrow Y\in F_Q$ if it is a Serre fibration, i.e. if it has the right lifting property with regard to the inclusions $\iota^q_0:I^q\rightarrow I^q\times I$ for $q\in\mathds N$.
    
    \item $f:X\rightarrow Y\in C_Q$ if $f\in{}^\boxslash (F_Q\cap W_Q)=\{\text{retracts of relative CW-complexes}\}$.
\end{itemize}

All objects in this model category are fibrant, and cofibrant objects are retracts of $CW$-complexes. In \texttt{Top} the model categorial notion of homotopy equivalence coincides with the classical one.

\begin{defin}
    A \textit{closed monoidal model category} is a category $\mathcal T$ equipped with both a closed monoidal structure and a model structure such that
    \begin{enumerate}
        \item Pushout-product axiom: for every pair of cofibrations $\iota:A\rightarrow X$ and $\iota':A'\rightarrow X'$ the canonical morphism
        $$
            (\iota\otimes 1_{X'},1_X\otimes \iota'):(A \otimes X') \coprod_{A \otimes A'} (X \otimes A')\rightarrow X \otimes X'
        $$
        is itself a cofibration, trivial if either $f$ or $f'$ is trivial;
        
        \item Unit axiom: The unit $\mathds 1$ is cofibrant.
    \end{enumerate}
    If the monoidal structure on $\mathcal T$ is closed symmetric (cartesian) then $\mathcal T$ is a \textit{closed symmetric (cartesian) model category}.
\end{defin}

The category $\texttt{Top}$ with the Quillen model structure and the closed cartesian structures presented here forms a closed cartesian model category \cite[4.2]{HoM07}.

Functors between model categories compatible with the model structures in the following way induce functors between their homotopy categories.

\begin{defin}
    An adjunction $(S\dashv\Lambda)$ is a \textit{Quillen adjunction} if the following equivalent conditions are satisfied:
    \begin{enumerate}
        \item $S$ preserves cofibrations and trivial cofibrations;
        
        \item $\Lambda$ preserves fibrations and trivial fibrations;
        
        \item $S$ preserves cofibrations and $\Lambda$ preserves fibrations;
        
        \item $S$ preserves trivial cofibrations and $\Lambda$ preserves trivial fibrations.
    \end{enumerate}
    
    An adjunction is a \textit{weak Quillen adjunction} if $S$ preserves cofibrant objects and weak equivalences between cofibrant objects and if $\Lambda$ preserves fibrant objects and weak equivalences between fibrant objects.
\end{defin}

Note that by Ken Brown's lemma \cite[Lem. 1.1.8]{HoM07} Quillen adjunctions are weak Quillen adjunctions.

\begin{defin}
    Let $S:\mathcal T\rightarrow \mathcal A$ be a functor that preserves cofibrant objects and weak equivalences between them and $\Lambda:\mathcal A\rightarrow T$ a functor that preserves fibrant objects and weak equivalences between fibrant objects.
    \begin{enumerate}
        \item The \textit{left derived functor} $\mathbb L S:\mathcal Ho\mathcal T\rightarrow \mathcal Ho\mathcal A$ is defined on objects as $\mathbb L S(X):=\text{Fib }S(X)$ and on morphisms as $\mathbb L S([f]):=[\text{Fib }S(f)]$.
        
        \item The \textit{right derived functor} $\mathbb R \Lambda:\mathcal Ho\mathcal A\rightarrow \mathcal Ho\mathcal T$ is defined on objects as $\mathbb R \Lambda(X):=\text{Cof }\Lambda(X)$ and on morphisms as $\mathbb R \Lambda([f]):=[\text{Cof }\Lambda(f)]$.
    \end{enumerate}
\end{defin}

If the two functors in the above definition are adjoint they induce an adjunction between the homotopy categories.

A useful tool for constructing new model categories out of known ones is by transfer through adjunctions. This requires that the known model structure is defined by a set of morphisms in the following sense.

\begin{defin}
    A model structure on a category $\mathcal T$  is \textit{cofibrantly generated} if there exists two sets of morphisms $I$ and $J$, called the \textit{generating cofibrations} and \textit{generating trivial cofibrations} respectively, such that
    \begin{enumerate}
        \item the domains of the morphisms in $I$ (resp. $J$) are small relative to transfinite composition of pushouts of coproducts of elements in $I$ (resp. $J$);
        
        \item $F=J^\boxslash$;
        
        \item $C={}^\boxslash(I^\boxslash)$.
    \end{enumerate}
\end{defin}

This condition allows the use of the small object argument in order to construct weak factorization systems. The Quillen model structure on $\texttt{Top}$ is cofibrantly generated. The set of generating cofibrations are the inclusions of the boundaries $\partial I^q\rightarrow I^{q}$ for $q\in\mathds N$ and the generating trivial cofibrations are the inclusions $I^q\rightarrow I^{q}\times I$ for $q\in\mathds N$.

\begin{prop}
    Let $\mathcal T$ and $\mathcal A$ be bicomplete categories, $(W,C,F)$ a cofibrantly generated model structure on $\mathcal T$ with a set of generating cofibrations $I$ and a set of generating trivial cofibrations $J$, and $(S\dashv\Lambda)$ be an adjunction between $\mathcal T$ and $\mathcal A$. Then sufficient conditions for $(\Lambda^{-1}(W),{}^\boxslash \Lambda^{-1}(W\cap F),\Lambda^{-1}(F))$ to be a cofibrantly generated model structure with generating cofibrations $S(I)$ and trivial cofibrations $S(J)$ is
    \begin{enumerate}
        \item $S$ preserves small objects;
        
        \item Any sequential colimit of pushouts of morphisms in $S(J)$ is contained in $\Lambda^{-1}(W)$.
    \end{enumerate}
    
    When these conditions hold, $(S\dashv\Lambda)$ is a Quillen adjunction.
\end{prop}

The $(-_+\dashv U)$-adjunction between $\texttt{Top}$ and $\texttt{Top}_\ast$ induces a cofibrantly generated closed monoidal model structure on $\texttt{Top}_\ast$.

\subsection{Loop space functors and sequential spectra}

Let $N\in\mathds N$ and $\mathds S^N$ denote the pointed space $I^N/_{\partial I^N}$ with base point the identified equivalence class $\partial I^N$.

\begin{defin}
    The \textit{$N$-fold loop space functor} is
    \begin{align*}
    \Omega^N:\texttt{Top}_\ast&\rightarrow\texttt{Top}_\ast\\
    X&\mapsto X^{\mathds S^N}
\end{align*}
\end{defin}

These functors admit Quillen left adjoints.

\begin{defin}
    The \textit{$N$-fold suspension functor} is
    \begin{align*}
    \Sigma^N:\texttt{Top}_\ast&\rightarrow\texttt{Top}_\ast\\
    X&\mapsto X\wedge \mathds S^{N}
\end{align*}
\end{defin}

\begin{prop}\label{NSigOmgAdj}
    The pair $(\Sigma^N\dashv \Omega^N)$ forms a Quillen adjunction.
\end{prop}

\textbf{Proof:} The unit of the adjunction is
\begin{align*}
    \eta^N:X&\rightarrow \Omega^N\Sigma^{N}X\\
    x&\mapsto \left(t\mapsto [x,t]\right),
\end{align*}
and the counit is
\begin{align*}
    \epsilon^N:\Sigma^{N}\Omega^N X &\rightarrow X\\
    [\gamma,t]&\mapsto \gamma(t)
\end{align*}
It is straightforward to verify that the counit-unit equations hold.

The pushout-product axiom implies that tensoring with a cofibrant objects preserves cofibrations and trivial cofibrations. Since $\mathds S^N$ is cofibrant and $\Sigma^N=-\wedge \mathds S^N$ it preserves cofibrations and trivial cofibrations.$\blacksquare$

In order to define the notion of infinite loop space we require the category of sequential spectra equipped with the stable model structure \cite{BF78} \cite{Sc97}.

\begin{defin}
    The category $\texttt{Sp}$ of \textit{sequential spectra} has as objects sequences of pointed topological spaces $X_\bullet\in\texttt{Top}_\ast^{\mathds N}$ equipped with structural maps $\sigma_\bullet:X_\bullet\wedge \mathds S^1\rightarrow X_{\bullet+1}$.
    
    \textit{Morphisms of sequential spectra} are sequences of pointed maps $f_\bullet:X_\bullet\rightarrow Y_\bullet$ such that the squares
    $$
    \xymatrix{
    X_\bullet\wedge\mathds S^1\ar[d]_{\sigma_\bullet}\ar[r]^{f_\bullet\wedge 1_{\mathds S_1}}&Y_\bullet\wedge\mathds S^1\ar[d]_{\sigma'_\bullet}\\
    X_{\bullet+1}\ar[r]_{f_{\bullet+1}}&Y_{\bullet+1}
    }
    $$
\end{defin}
are commutative.

The stable homotopy groups of sequential spectra are defined as
$$
    \pi_q(X_\bullet)=\varinjlim\pi_{\bullet+q}(X_\bullet)
$$

The stable model structure $(W_S,C_S,F_S)$ on $\texttt{Sp}$ is given by the following classes of maps:
\begin{itemize}
    \item $f_\bullet:X_\bullet\rightarrow Y_\bullet\in W_S$ if $f_{\bullet\ast}:\pi_q(X_\bullet)\rightarrow \pi_q(Y_\bullet)$ is an isomorphism for all $q\in\mathds Z$.
    
    \item $f_\bullet:X_\bullet\rightarrow Y_\bullet\in F_S$ if if each $f_q:X_q\rightarrow Y_q$ is in $F_Q$, and the maps
    $$
        (\Omega^1(\sigma^X_q)\eta^1_{X_q},f_q):X_q\rightarrow \Omega^1 X_{q+1}\times_{\Omega^1 Y_{q+1}}Y_q
    $$
    are in $W_Q$.
    
    \item $f_\bullet:X_\bullet\rightarrow Y_\bullet\in C_S$ if $f_0:X_0\rightarrow Y_0$ and each map
$$
[f_{q+1},\sigma_q^Y]:X_{q+1}\cup_{X_q\wedge \mathds S^1}Y_q\wedge \mathds S^1\rightarrow Y_{q+1}
$$
are in $C_Q$.
\end{itemize}

The sequential spectra of interest for stable homotopy theory are the $\Omega$-spectra, which are the ones such that $\Omega^1(\sigma_\bullet)\eta^1_{X_\bullet}:X_\bullet\rightarrow \Omega^1 X_{\bullet+1}$ are weak equivalences. These are the fibrant objects of the stable model structure on sequential spectra. The cofibrant objects are sequential spectra where $X_0$ is a retract of a CW-complex and each structural map is a cofibration in $\texttt{Top}_\ast$.

\begin{defin}
    A spectrum $X_\bullet$ is \textit{$m$-connective} if $\pi_q(X_\bullet)$ is trivial for all $q\leq m$.
    
    We will denote the category of $m$-connective spectra as $\texttt{Sp}_m$.
\end{defin}

\begin{defin}
    The \textit{$\infty$-fold loop space functor} is
    \begin{align*}
    \Omega^\infty:\texttt{Sp}&\rightarrow\texttt{Top}_\ast\\
    X_\bullet&\mapsto \varinjlim\Omega^{\bullet} X_\bullet
\end{align*}
\end{defin}

This functor admits a Quillen left adjoint.

\begin{defin}
    The \textit{$\infty$-fold suspension functor} is
    \begin{align*}
    \Sigma^\infty:\texttt{Top}_\ast&\rightarrow\texttt{Sp}\\
    X&\mapsto X\wedge \mathds S^{\bullet}
\end{align*}
\end{defin}

Proposition \ref{NSigOmgAdj} implies the following.

\begin{cor}
    The pair $(\Sigma^\infty\dashv \Omega^\infty)$ forms a Quillen adjunction.
\end{cor}

\subsection{Relative spaces, relative loop spaces and Relative spectra}

Let $\mathcal T$ be a model category with model structure $(W,F,C)$. Then there is the projective model structure $(W_{\text{proj}},C_{\text{proj}},F_{\text{proj}})$ on the category of morphisms $\mathcal T^{\shortrightarrow}$ \cite[Prop. A.2.8.2]{Lu09}, whose objects are morphisms $\iota:A\rightarrow X$ of $\mathcal T$ and morphisms are commutative squares of morphisms
$$
    \xymatrix{
        A\ar[r]^f\ar[d]_\iota&B\ar[d]^{\iota'}\\
        X\ar[r]_{g}&Y
    }
$$
called the projective model structure. In this model structure
\begin{enumerate}
    \item $(f,g):\iota\rightarrow \iota'\in W_{\text{proj}}$ if $f,g\in W$;
    
    \item $(f,g):\iota\rightarrow \iota'\in F_{\text{proj}}$ if $f,g\in F$;
    
    \item $(f,g):\iota\rightarrow \iota'\in C_{\text{proj}}$ if $f,[g,\iota']\in C$, where $[g,\iota']$ is the unique morphism
    $$
        [g,\iota']:X\cup_{A} B\rightarrow Y
    $$
    given by the pushout property.
\end{enumerate}

Fibrant objects in this category are morphisms between fibrant objects, and cofibrant objects are cofibrations between cofibrant objects.

Let $\texttt{Top}_\ast^{\shortrightarrow}$ be the category of pointed maps with the projective model structure. We call it the category of \textit{relative spaces}. All objects in this category are fibrant, and the cofibrant objects are retracts of inclusions of pointed CW-pairs.

\begin{defin}
    A relative space $\iota:A\rightarrow X\in\texttt{Top}_\ast^{\shortrightarrow}$ is \textit{$m$-connected} if $X$ is $m$-connected and $A$ is $m-1$-connected.
    
    We will denote the category of $m$-connected relative spaces as $\texttt{Top}^\shortrightarrow_m$.
\end{defin}

Let $N\in\mathds N$ and $I$ be the pointed space $([0,1],1)$.

\begin{defin} The \textit{$N$-fold  relative loop space functor} is
    \begin{align*}
        \Omega^N_\text{rel}:\texttt{Top}_\ast^{\shortrightarrow}&\rightarrow \texttt{Top}_\ast\\
        \iota:A\rightarrow X&\mapsto (A\times_X X^{ I})^{\mathds S^{N-1}}
    \end{align*}
    where $A\times_X X^{ I}$ is the pullback:
    $$
    \xymatrix{
    A\times_X X^{ I}\ar[r]\ar[d]&X^I\ar[d]^{\text{ev}_0}\\
    A\ar[r]_\iota&X
    }
$$
i.e. the space $\{(a,\beta)\in A\times X^I\ |\ \beta(0)=a,\ \beta(1)=\ast_X\}$. This construction is known as the homotopy fiber of $\iota$.
\end{defin}

There is an obvious functor $\text{incl}:\texttt{Top}_\ast\rightarrow \texttt{Top}^\shortrightarrow_\ast$, and $\Omega^N_\text{rel}\text{ incl}=\Omega^N$. If $\iota:A\rightarrow X$ is $m$-connected then $\Omega_\text{rel}\iota$ is $m-1$-connected.

\begin{defin} The \textit{$N$-fold 2-loop space functor} is
    \begin{align*}
        \Omega^N_2:\texttt{Top}^\shortrightarrow_{\ast}&\rightarrow \texttt{Top}^2_{\ast}\\
        \iota:A\rightarrow X&\mapsto (\Omega^N_\text{rel}\iota,\Omega^N X).
    \end{align*}
\end{defin}

This functor admits a weakly Quillen left adjoint. As we will later be working with 2-operads with $2$ denoting the ordered set $\{o< c\}$, and these colored operads will act on the category of pairs of spaces we use the notation $(X_o,X_c)$ for objects of $\texttt{Top}_\ast^2$.

\begin{defin}
    The \textit{$N$-fold 2-suspension functor} is
\begin{align*}
    \Sigma_2^N:\texttt{Top}_\ast^2&\rightarrow\texttt{Top}^\shortrightarrow_{\ast}\\
    (X_o,X_c)&\mapsto \left(\begin{aligned}X_o\wedge\mathds S^{N-1}&\rightarrow (X_o\wedge I\wedge\mathds S^{N-1})\vee  (X_c\wedge\mathds S^N)\\
    [x_o,s]&\mapsto[x_o,0,s]\end{aligned}\right)
\end{align*}
\end{defin}

Note that for $f:X\rightarrow Y$ a cofibration in $\texttt{Top}_\ast$ it is not generally the case that $Y\cup_XX\wedge I\rightarrow Y\wedge I$ is a cofibration. This implies that $\Sigma^N_2$ doesn't preserve cofibrations in general. Also $\Omega^N_\text{rel}$  doesn't preserve fibrations, and therefore neither does $\Omega^N_2$.

\begin{prop}\label{2NSigOmgAdj} The pairs $(\Sigma_2^N\dashv \Omega_2^N)$ form weak Quillen adjunctions.
\end{prop}

\textbf{Proof:} The unit of the adjunction is
\begin{align*}
    \eta^N_2:(X_o,X_c)&\rightarrow \Omega^N_2\Sigma^N_2(X_o,X_c)\\
    x_o&\mapsto \left(s\mapsto ([x_o,s],s'\mapsto [x_o,s',s])\right)\\
    x_c&\mapsto \left(t\mapsto [x_c,t]\right)
\end{align*}
and the counit is
\begin{align*}
    \epsilon^N_2:\Sigma^N_2\Omega^N_2\iota &\rightarrow \iota\\
    [(\alpha,\beta),s]&\mapsto \alpha(s)\\
    [(\alpha,\beta),s',s]&\mapsto\beta(s)(s')\\
    [\gamma,t]&\mapsto\gamma(t).
\end{align*}
It is straightforward to verify that the counit-unit equations hold.

Since tensoring with cofibrant objects preserves cofibrations, and therefore preserves cofibrant objects, and coproducts of cofibrant objects are cofibrant, the image of $\Sigma^N_2$ on cofibrant objects is a map between cofibrant objects. Also the inclusion of a cofibrant objects on the base of its cone is a cofibration. Therefore $\Sigma^N_2$ preserves cofibrant objects. Both tensoring with a cofibrant object and taking coproducts preserves weak equivalences between cofibrant objects, therefore $\Sigma^N_2$ preserves weak equivalences between cofibrant objects. 

Clearly $\Omega^N_2$ preserves fibrant objects since all objects in $\texttt{Top}_\ast^2$ are fibrant. Since the $\Omega^N$ are right Quillen functors they preserve weak equivalences. Now note that for any relative space $\iota:A\rightarrow X$ we get a natural exact sequence of pointed spaces
$$
    \Omega^NA\rightarrow \Omega^NX\rightarrow \Omega^N_\text{rel}\iota\rightarrow \Omega^{N-1}A\rightarrow \Omega^{N-1}X
$$
which induce an exact sequence on homotopy groups \cite[Ch 8.6]{Ma99} and therefore by the five lemma $\Omega^N_\text{rel}$ preserves weak equivalences.$\blacksquare$

We call $\texttt{Sp}^\shortrightarrow$ with the projective model structure the category of relative spectra.

\begin{defin}
    A relative spectrum $\iota_\bullet:A_\bullet\rightarrow X_\bullet$ is \textit{$m$-connective} if $X_\bullet$ is $m$-connective and $A_\bullet$ is $m-1$-connective.
    
    We will denote the category of $m$-connective relative spectra as $\texttt{Sp}^\shortrightarrow_m$.
\end{defin}

\begin{defin}
    The \textit{$\infty$-fold 2-loop space functor} is
    \begin{align*}
        \Omega^\infty_2:\texttt{Sp}^\shortrightarrow&\rightarrow \texttt{Top}^2_{\ast}\\
        \iota_\bullet:A_\bullet\rightarrow X_\bullet&\mapsto(\varinjlim \Omega^{\bullet+1}_\text{rel}\iota_\bullet,\varinjlim\Omega^{\bullet+1} X_\bullet).
    \end{align*}
\end{defin}

\begin{defin}
    The \textit{$\infty$-fold 2-suspension functor} is
\begin{align*}
    \Sigma_2^\infty:\texttt{Top}_\ast^2&\rightarrow\texttt{Sp}^\shortrightarrow\\
    (X_o,X_c)&\mapsto \left(X_o\wedge\{0\}\wedge\mathds S^{\bullet}\rightarrow (X_o\wedge I\wedge\mathds S^{\bullet})\vee  (X_c\wedge\mathds S^{\bullet+1})\right)
\end{align*}
\end{defin}

Proposition \ref{2NSigOmgAdj} implies the following.

\begin{cor}
    The pair $(\Sigma_2^\infty\dashv \Omega_2^\infty)$ forms a weak Quillen adjunction.
\end{cor}

\section{Monads, operads and 2-operads}

A monad is an endofunctor equipped with certain natural transformations that allows us to interpret its images as free objects of an algebraic theory. Each monad defines algebras which are the objects of its algebraic theory. Operads are structures that allow us to define monads in a category using objects of the category itself. 2-operads allow us to define algebraic structures on pairs of objects, one of which acts on the other.

The $(\Omega^N\dashv\Sigma^N)$ and $(\Omega^N_2\dashv\Sigma^N_2)$ adjunctions give $\Omega^N\Sigma^N$ and $\Omega^N_2\Sigma^N_2$ the structure of monads in $\texttt{Top}_\ast$ and $\texttt{Top}_\ast^2$ respectively. As will be seen in the next section the little cubes operads $\mathcal C^N$ and the Swiss-cheese 2-operads $\mathcal{SC}^N$ have monads associated to them that are closely related to the monads of these adjunctions.

In 3.1  the definitions of monads, algebras over monads and of monad functors are given. In 3.2 the two sided bar construction and the geometric realization functor are defined. These constructions define all spaces and maps needed to prove the relative recognition principle. In 3.3 operads and their associated monads are presented. In 3.4 the little cubes operads $\mathcal C^N$ are defined and the $\mathcal C^N$ actions on $N$-fold loop spaces are shown. In 3.5 2-operads and their associated monads are defined. In 3.6 the Swiss-cheese 2-operads $\mathcal{SC}^N$ are defined and the $\mathcal{SC}^N$-actions on $N$-fold 2-loop spaces are shown. In 3.7 we present a review of the theory of model structures on 2-operads and their algebras.

\subsection{Monads}

In \cite{Ma72} the language of monads is central in the proof of the recognition principle.

\begin{defin}
    Let $\mathcal T$ be a category. The \textit{endofunctor monoidal category of $\mathcal T$} is the category $\texttt{Cat}(\mathcal T,\mathcal T)$, the category of endofunctors on $\mathcal T$, equipped with composition as the tensor product and the identity functor as the identity.
\end{defin}

\begin{defin}
    A \textit{monad in $\mathcal T$} is a monoid $(C,\mu,\eta)$ in $\text{End}(\mathcal T)$, i.e. it is an endofunctor $C:\mathcal T\rightarrow \mathcal T$ equipped with natural transformations $\mu:CC\rightarrow C$ and $\eta:Id_{\mathcal T}\rightarrow C$ such that the diagrams below commute:
    
    $$
        \xymatrix{
            CCC\ar[r]^{1_C \mu}\ar[d]_{\mu 1_C}&CC\ar[d]^\mu&&&C\ar[r]^{1_{C}\eta}\ar[dr]_{1_C}&CC\ar[d]|\mu&C\ar[l]_{\eta{1_C}}\ar[dl]^{1_C}\\
            C C\ar[r]_\mu&C&&&&C&
        }
    $$
\end{defin}

\begin{defin}
    A \textit{$C$-algebra over a monad $C$ in $\mathcal T$} is an object $X\in \mathcal T$ equipped with a morphism $\xi:CX \rightarrow X$ such that the diagrams bellow commute:
    $$
        \xymatrix{
            C CX\ar[r]^{C\xi}\ar[d]_{\mu_X}&CX\ar[d]^\xi&&&X\ar[r]^{\eta_X}\ar[dr]_{1_X}&CX\ar[d]^\xi\\
            CX\ar[r]_\xi&X&&&&X
        }
    $$
    
    A morphism of $C$-algebras is a morphism between the objects that commute with the structural morphisms. We refer to them as C-maps.
    
    We denote the category of $C$-algebras as $C[\mathcal T]$. When $\mathcal T=\texttt{Top}$ we refer to $C$-algebras as $C$-spaces.
\end{defin}

Note that for any $X\in\mathcal T$, $CX\in C[\mathcal T]$. Therefore $C$ defines a functor $C:\mathcal T\rightarrow C[\mathcal T]$. This functor is left adjoint to the forgetful functor.

\begin{defin}
    A \textit{$C$-functor in $\mathcal A$ of a monad $C$ in $\mathcal T$} is a functor $F:\mathcal T\rightarrow\mathcal A$ equipped with a natural transformation $\lambda: F\circ C\rightarrow F$ such that the diagrams bellow commute:
    $$
        \xymatrix{
            F C C\ar[r]^{\lambda_C}\ar[d]_{F\circ \mu }& F C\ar[d]^\lambda&&&F\ar[r]^{F\circ \eta }\ar[dr]_{1_F}&F C\ar[d]^\lambda\\
            F C\ar[r]_\lambda&F&&&&F
        }
    $$
    
    A morphism of $C$-functors in $\mathcal A$ is a natural transformation between the functors that commute with the structural natural transformations.
\end{defin}

The examples we will need fall into the following cases:
\begin{enumerate}
    \item If $(C,\mu,\nu)$ is a monad then $(C,\mu)$ is a $C$-functor.
    
    \item If $\phi: (C,\mu,\nu)\rightarrow (C',\mu',\nu')$ is a morphism of monads and $(F,\lambda)$ is a $C'$-functor then $(F,\lambda\circ F\phi)$ is a $C$-functor. In particular, by the previous example, $(C',\mu'\circ C'\phi)$ is a $C$-functor. We also have that $\phi:(C,\mu)\rightarrow (C',\mu' C'(\phi))$ is a morphism of $C$-functors.
    
    \item Let $(S\dashv \Lambda):\mathcal T\rightleftharpoons\mathcal A$ be an adjunction with unit $\eta:Id_{\mathcal T}\rightarrow \Lambda S$ and counit $\epsilon:S\Lambda\rightarrow Id_{\mathcal A}$. Then $(\Lambda S,\Lambda\epsilon_S,\eta)$ is a monad in $\mathcal T$ and $(S,\epsilon_S)$ is a $\Lambda S$-functor.
    
    \item If $\alpha:(C,\mu,\nu)\rightarrow (\Lambda S,\Lambda\epsilon_S,\eta)$ is a morphism of monads, with $\Lambda S$ as in the previous example, then $S$ is a $C$-functor in $\mathcal A$ and $\alpha:(C,\mu)\rightarrow (\Lambda S,\Lambda(\epsilon_{S}) S(\alpha))$ is a morphism of $C$-functors in $C[\mathcal T]$.
\end{enumerate}

\subsection{Bar construction}

Remember that the simplicial category $\Delta$ has as objects the linearly ordered sets $\langle q\rangle:=\{0<1<\dots<q\}$ and morphisms are order preserving functions between these ordered sets. The morphisms in this category are generated by morphisms $\partial_i:\langle q-1\rangle\rightarrow \langle q \rangle$ and $s_i:\langle q+1\rangle\rightarrow \langle q\rangle$ for each $i\in\langle q\rangle$. A simplicial object in $\mathcal T$ is a functor $X_\bullet:\Delta^\text{op}\rightarrow \mathcal T$, and a cosimplicial object in $\mathcal T$ is a functor $X^\bullet:\Delta\rightarrow \mathcal T$.

The category of simplicial objects on a model category $\mathcal T$ admits a model structure called the Reedy model structure (see \cite[5.2]{HoM07} for details).

The two sided bar construction described here was introduced in \cite{Ma72} and admits a variety of applications.

\begin{defin}
    The \textit{two sided bar construction} is the category $B(\mathcal T,\mathcal A)$ and the functor $B_\bullet : B(\mathcal T,\mathcal A) \rightarrow \mathcal A^{\Delta^\text{op}}$ defined as follows. The objects of $B(\mathcal T,\mathcal A)$ are triples $((F,\lambda),(C,\mu,\nu),(X,\xi))$, abbreviated $(F,C,X)$, where $C$ is a monad in $\mathcal T$, $F$ is a $C$-functor in $\mathcal A$ and $X$ is a $C$-algebra, and morphisms are triples $(\alpha,\phi,f):(F,C,X)\rightarrow(F',C',X')$ where $\phi:C\rightarrow C'$ is a monad morphism, $f:X\rightarrow X'$ is a $C$-algebra morphism and $\alpha:F\rightarrow F'$ is a morphism of $C$-functors. 
    
    The functor $B_\bullet(F,C,X)$ is
    \begin{align*}
        B_\bullet:B(\mathcal T,\mathcal A)&\rightarrow \mathcal A^{\Delta^\text{op}}\\
        (F,C,X)&\mapsto\left(FC^\bullet X,\ \begin{matrix*}[l]
        \partial_i=\begin{cases}
            \lambda_{C^{\bullet-1}},&i=0\\
            F C^{i-1}\mu_{C^{\bullet-i}},&0<i<\bullet\\
            F C^{\bullet-1}\rho,&i=\bullet
        \end{cases}\\
        s_i=FC^{ i}\eta_{C^{ \bullet-i+1}}, 0\leq i\leq \bullet\end{matrix*}\right)
    \end{align*}
\end{defin}

Given a small category $\mathcal C$ and a functor $F:\mathcal C^{\text{op}}\times\mathcal C\rightarrow \mathcal T$ we can define the coend of $F$ as the following coequilizer:
$$
    \int^{\mathcal C}F(X,X):= \text{Coeq}\left(\coprod_{X\rightarrow X'}F(X',X)\rightrightarrows\coprod_{X\in \mathcal C}F(X,X)\right)
$$

Let $\Delta^\bullet\in\texttt{Top}^\Delta$ be the simplex cosimplicial space where
$$
\Delta^\bullet := 
  \left\{
    u \in \mathds{R}^{\bullet}
    \left|
    \sum_{i \in \bullet } x^i = 1 ;
    \forall i\in \bullet:x^i \geq 0 
  \right.\right\}.
  $$

\begin{defin}
    The geometric realization functor is
    \begin{align*}
        |-|:\texttt{Top}^{\Delta^{\text{op}}}&\rightarrow\texttt{Top}\\
        X_\bullet&\mapsto \int^{\Delta}X_q\times \Delta^q 
    \end{align*}
\end{defin}

This functor induces a geometric realization functor on $\texttt{Top}_\ast$, $\texttt{Top}_\ast^2$, $\texttt{Top}_\ast^{\shortrightarrow}$, $\texttt{Sp}$ and $\texttt{Sp}^{\shortrightarrow}$, and we denote by $B(F,C,X)$ the geometric realization $|B_\bullet(F,C,X)|$ when $F$ is a $C$-functor on any of these categories. From \cite[9.2, 11.8]{Ma72} we have that any map $f:Y\rightarrow FX$ determines a map $\tau(f):Y\rightarrow B(F,C,X)$ and any map $g:FX\rightarrow Y$ such that $g\partial_0=g\partial_1:FCX\rightarrow Y$ determines a map $\varepsilon(g):B(F,C,X)\rightarrow Y$.

We will need the following definitions.

\begin{defin}
     For $X$ a space and $A$ a closed subspace of $X$ we say that $(A,X)$ is a \textit{neighborhood deformation retract pair}, or \textit{NDR-pair}, if there are maps $u:X\rightarrow I$ and $H:X\times I\rightarrow X$ such that $u^{-1}(0)\subset A$, $H(x,0)=x$ for all $x\in X$, $H(a,t)=a$ for all $t\in I$ and $a\in A$, and $H(x,1)\in A$ for all $x\in u^{-1}([0,1))$.
\end{defin}

\begin{defin}
     A pointed space $X\in\texttt{Top}_\ast$ is well-pointed if $(\ast, X)$ is an NDR-pair.
\end{defin}

Note that any cofibrant pointed space is well-pointed, so any space can be replaced by a weakly equivalent well-pointed space.

\begin{defin}
    A simplicial space $X_\bullet \in\texttt{Top}^{\Delta^{\text{op}}}$ is \textit{proper} if $\text{colim}_{s_i}X_{q-1}\rightarrow X_q$ is an NDR-pair for all $q\in\mathds N$, where $\text{colim}_{s_i}X_{q-1}$ is the union of the images of the degeneracies $s_i$.
\end{defin}

\subsection{Operads}

Operads were introduced in \cite{Ma72} in order to study loop spaces, though the operad of associahedras were already implicitly described in \cite{StJ63}.

\begin{defin}
    Let $\mathbb S_\text{inj}$ be the category whose objects are (possibly empty) finite sets $S$ and whose morphisms are injective functions. We will denote by $\text{sk}\mathbb S_\text{inj}$ the skeleton of $\mathbb S_\text{inj}$ whose object are $\underline{n}:=\{1,\dots,n\}$ for $n\in\mathds N$, where $\underline{0}:=\emptyset$. 

    Let $\mathbb S$ be the subcategory of $\mathbb S_\text{inj}$ with the same objects and bijective functions as morphisms, and $\text{sk}\mathbb S$ its skeleton.
\end{defin}

For every $S\in \mathbb S$ the bijections $\mathbb S_S:=\mathbb S(S,S)$ have a group structure.

\begin{defin}
    Let $(\mathcal T,\otimes,\mathds 1)$ be a bicomplete symmetric monoidal category. An \textit{$\mathbb S$-object in $\mathcal T$} is a functor $\mathcal O:\text{$\mathbb S^\text{op}$}\rightarrow\mathcal T$ such that $\mathcal O$ is reduced, in the sense that $\mathcal O(\underline 0)$ is the terminal object $\ast$.
\end{defin}

We note that the terminology for a "reduced" $\mathbb S$-object is inconsistent in the literature. We follow the terminology used in \cite{Ma09}, which is consistent with the definition for operads in \cite{Ma72}. This differs from the notion of "reduced" $\mathbb S$-object in \cite{BM03} where "reduced" is used for $\mathbb S$-objects where $\mathcal O(\underline 0)$ is the unit $\mathds 1$ of the monoidal structure. In cartesian categories, where the unit is the terminal object, the two notions coincide.

For an $\mathbb S$-object in $\mathcal T$ and $S\in\mathbb S$ we get an action of $\mathbb S_{S}$ on $\mathcal O(S)$. An $\mathbb S$-object in $\texttt{Top}$ is referred to as an $\mathbb S$-space.

\begin{defin}
     An $\mathbb S$-space $\mathcal O$ is $\mathbb S$-free if the $\mathbb S_{S}$ actions on $\mathcal O(S)$ are free.
\end{defin}

We can use the language of trees to define a monad in the category of $\mathbb S$-objects.

\begin{defin}
    A \textit{graph} $T$ is a 1-dimensional simplicial complex, not necessarily compact (that means an edge doesn't need to end at a vertex, these edges are called \textit{external}. Edges that end at vertices at both sides are called \textit{internal}). The set of vertices and edges of $T$ are denoted by $V_T$ and $E_T$ respectively.
    
    An \textit{orientation of a graph} is an orientation on each edge. For $v\in V_T$ we denote by $\text{in}(v)$ the set of incoming edges and by $\text{out}(v)$ the set of outgoing edges incident to $v$. We also denote by $\text{in}(T)$ and by $\text{out}(T)$ the sets of incoming and outgoing external edges, respectively. We note that $\text{in}(v)$ can be empty.
    
    A \textit{tree} is an oriented, contractible, finite graph such that $\text{out}(T)$ has one element denoted $e_0$ and for all $v\in V_T$ the set $\text{out}(v)$ has one element. We call the elements of $\text{in}(T)$ the \textit{leaves} of the tree. Notice that there are partial orders on both the vertices and the edges such that $x\leq x'$ if $x$ is in the unique path from $x'$ to $e_0$ following the orientation of the tree. Define $\text{pre}(v):=\max\{v'\in V_T\ |\ v'\leq v;\ v\neq v' \}$.
    
    An \textit{ordering on a tree} is a linear order of $V_T$.
    
    A \textit{tree on a finite set $S$} is a tree $T$ equipped with a bijection between $\text{in}(T)$ and $S$, called the \textit{labeling of the leafs}.
    
    We denote by $\mathbb T(\mathbb S)$ the category whose objects are isomorphism classes of ordered trees on finite sets and whose morphisms are (not necessarily order preserving) isomorphisms of trees compatible with the labeling of the leafs. We also define $\mathbb T(\text{sk}\mathbb S)$ the subcategory of trees on the sets in $\text{sk}\mathbb S$.
    
    We refer to the unique trees on $S$ with only one vertex as the $S$-corolla. Note that there is a unique tree with no vertices and only one edge.
\end{defin}

There is a grafting operation on trees. Let $T_0\in \mathbb T(S)$ and $T_v\in \mathbb T(\text{in}(v))$ for $v\in V_{T_0}$, we define the tree $T_0[\{T_v\}_{v\in V_{T_0}}]\in\mathbb T(S)$ as the one obtained from the union of the $T_v$ by identifying the root of each $T_v$ with the leaf of $T_{\text{pre}(v)}$ corresponding to $\text{out}(v)$. If $T_0$ and each $T_v$ are ordered then $T_0[\{T_v\}_{v\in V_{T_0}}]$ inherits an order.

\begin{defin}
    For each $\mathbb S$-object $\mathcal O$ define the $\mathbb S$-object
    \begin{align*}
        \mathbb T\mathcal O: \mathbb{S}^{\text{op}}&\rightarrow \mathcal T\\
        S&\mapsto \text{Coeq}\left(\coprod_{T'\rightarrow T}\bigotimes_{v\in V_{T'}}\mathcal O(\text{in}(v))\rightrightarrows\coprod_{T\in \mathbb T(S)}\bigotimes_{v\in V_T}\mathcal O(\text{in}(v))\right)
    \end{align*}
    We need the trees to be ordered in order to define the tensor product, but we don't want the specific ordering to matter, and that is why we take the coequalizer. For $S=\underline 1$ the summand associated to the tree without vertices is $\mathds 1$.
    
    This defines an endofunctor on $\mathcal T^{\mathbb S^{\text{op}}}$ that is the underlying functor of a monad. The unit $\eta^{\mathbb T}_{\mathcal O}:\mathcal O\rightarrow \mathbb T\mathcal O$ is the natural inclusion of $\mathcal O(S)\rightarrow \mathbb T\mathcal O(S)$ in the summand indexed by the $S$-corolla. The product $\mu^{\mathbb T}_{\mathcal O}:\mathbb T\mathbb T\mathcal O\rightarrow\mathbb T\mathcal O$ is induced by the grafting operation on trees.
\end{defin}

\begin{defin}
    An operad $\mathcal O$ in $\mathcal T$ is a $\mathbb T$-algebra in $\mathcal T^{\mathbb S^{\text{op}}}$.
    
    The operad structure of $\mathcal O$ is completely defined by the morphisms
$$
\eta^{\mathcal O}:\mathds 1\rightarrow \mathcal O(\underline{1}),
$$
which is the composition of the inclusion of the $\mathds 1$ summand of $\mathbb T\mathcal O(1)$ associated with the tree without vertices with the $\mathbb T$-algebra structural morphism of $\mathcal O$, and, for each $P$ linearly ordered set and family of finite sets $S_p$ indexed by $P$, a morphism
$$
    \mu^{\mathcal O}_\nu:\mathcal O(P)\otimes \left(\bigotimes_{p\in P}\mathcal O(S_p)\right)\rightarrow\mathcal O(\cup_{p\in P}S_p),
$$
given by the $\mathbb T$-algebra structural morphism on the summand indexed by a tree obtained by identifying the roots of the $\nu^{-1}(p)$-corollas with the $p$ leaf of the $P$-corolla, with order given by the order on $P$. These morphisms satisfy associative, unit and equivariance laws.
\end{defin}

Any operad in $\mathcal T$ defines a monad in $\mathcal T_\ast$, the category of pointed objects of $\mathcal T$. If $\mathcal O$ is an operad we can extend the underlying functor on $\mathbb S^{\text{op}}$ to a functor on $\mathbb S_{\text{inj}}^{\text{op}}$. For an inclusion $\nu:S\hookrightarrow S'$ we can build a tree by adding a vertex with no incoming edges to the end of each leaf of the $S'$-corolla that is not in the image of $\nu$. This results in a tree on $S$ and the operad structure gives us a morphism $\mathcal O(S')\rightarrow \mathcal O(S)$. We refer to these morphisms as degeneracies of the operad $\mathcal O$.

Any $X\in \mathcal T_\ast$ defines a functor
\begin{align*}
    X^-:\mathbb S_{\text{inj}}&\rightarrow \mathcal T\\
    S&\mapsto X^S
\end{align*}
with the morphisms induced by the injective functions of $\mathbb S_\text{inj}$ exchanging the coordinate labels and inserting the base point $\ast_X$ on the coordinates that aren't contained in the image of the function.

We can then define for any operad $\mathcal O$ the endofunctor
\begin{align*}
    O:\mathcal T_\ast&\rightarrow \mathcal T_\ast\\
    X&\mapsto \int^{\text{sk}\mathbb S_{\text{inj}}} \mathcal O(\underline n)\times X^{\underline{n}}
\end{align*}
with base point given by $\mathcal O(\underline 0)\times X^{\underline{0}}=\ast$.

This endofunctor has a monad structure, with the unit given by $x\mapsto [\eta^{\mathcal O},x]$ and multiplication induced by the operad multiplication $\mu^{\mathcal O}$.

\begin{defin}
    For an operad $\mathcal O$ an \textit{$\mathcal O$-algebra} is an $O$-algebra for the associated monad $O$.
\end{defin}

If $\mathcal O$ is an operad in $\texttt{Top}$ we refer to $\mathcal O$-algebras as $\mathcal O$-spaces.

For every operad $\mathcal O$ and pointed object $(X,\ast_X)$ the pointed object $OX$ admits a natural filtration. Let $F^k\text{sk}\mathbb S_{\text{inj}}$ be the full subcategory of $\text{sk}\mathbb S_{\text{inj}}$ containing only the sets $\underline n$ with $n\leq k$. We can then define $F^k OX$ as the image of the natural inclusion of
$$
    \int^{F^k\text{sk}\mathbb S_{\text{inj}}}\mathcal O(\underline n)\times X^{\underline{n}}
$$
into $O(X)$.

\subsection{Little cubes operads}

The little $N$-cubes topological operads $\mathcal{C}^{N}$ were first introduced by Boardman and Vogt in the language of PROPS \cite{BV68}, and their operad structures is described by May in \cite{Ma72}. For $N<\infty$ a little $N$-cube is a linear embedding $d:I^N\rightarrow I^N$ with parallel axes. This means $d$ is of the form $d(t)=Mt+C$ for some diagonal $N\times N$ matrix $M$ with strictly positive entries and some $C\in I^N$. Define $\mathcal C^N(S)$ as the set of configurations of little $N$-cubes $d_S$ indexed by the elements of the set $S$ such that the images of the $d_a$ for $a\in S$ are pairwise disjoint. We can regard $d_S$ as an element of ${I^N}^{\coprod_{S}  I_a^{N}}$ and give $\mathcal C^N(S)$ the subspace topology. Note that $\mathcal C(\emptyset)=\{d_\emptyset\}$, where $d_\emptyset$ is the unique map $\emptyset \rightarrow I^N$.

The $\mathbb S$-space structural maps induced by the bijections of $S$ switch the indices of the $d_a$. The operad structure is given by composition of the little cubes. The degeneracies delete the little $N$-cubes indexed by the elements that are not in the image of the inclusions $S\rightarrow S'$.

We have natural operad inclusions
\begin{align*}
    \sigma^N_S:\mathcal{C}^N(S)&\rightarrow\mathcal{C}^{N+1}(S)\\
    \coprod_{a\in S} d_a&\mapsto\coprod_{a\in S} d_a \times1_I
\end{align*}
and we define $\mathcal{C}^\infty:=\varinjlim \mathcal{C}^N$, with the topology of the union.

Note that for all $S\in\mathbb S$ the spaces $\mathcal{C}^\infty(S)$ are contractible.

\begin{teor}
    The images of the $N$-fold loop space functors $\Omega^N$ are naturally $\mathcal C^N$-spaces.
\end{teor}

\textbf{Proof:} First assume $N<\infty$. Let $\theta^N:C^{N}\Omega^N\rightarrow \Omega^N$ be the natural transformation defined by
\begin{align*}
    \theta^N_X:C^{N}\Omega^N X&\rightarrow \Omega^NX\\
    \left[d_{\underline n},\gamma^{\underline n}\right]&\mapsto \left(t\mapsto\begin{cases}
    \gamma^i(d_i^{-1}(t)),&t\in d_i(I^N)\\
    \ast_c,&t\not\in d_{\underline n}\left(\coprod_{\underline n}I^N\right)
    \end{cases}\right)
\end{align*}
It is straightforward to verify that $\theta^N$ is compatible with the operad structural morphisms of $\mathcal{C}^N$.

The action of $\mathcal{C}^\infty$ on the image of $\Omega^\infty$ is induced from the finite cases. By definition an element $d_{\underline{n}}\in \mathcal{C}^\infty(\underline{n})$ belongs to some $\mathcal{C}^N(\underline{n})$ with $N<\infty$. We can then use the above construction for the finite cases to define the structural natural transformation $\theta^\infty$.$\blacksquare$

\subsection{2-operads}

The theory of colored operads is a generalization of operad theory where we allow operations on multiple objects. We will focus on the case with two colors. Let 2 be the ordered set $\{o<c\}$. Consider the functor from sets to pointed partially ordered sets
\begin{align*}
    -_\ast:\texttt{Set}&\rightarrow\texttt{PO-Set}_\ast\\
    S&\mapsto (\{S\}\coprod S, S)
\end{align*}
    where the partially ordered structure is given by $S<a$, $\forall a\in S$.

\begin{defin}
    A \textit{2-set} is a set $S$ equipped with an unbased order preserving function $\text{cor}_S:S_\ast\rightarrow 2$, which will be referred to as the coloring of $S$. A 2-set is usually simply denoted $S$ and the coloring $\text{cor}_S$ will be denoted simply as $\text{cor}$ when the 2-set $S$ is evident from the context.
    
    We also define $S_o:=\text{cor}^{-1}(o)\cap S$ the \textit{open part of $S$}, and $S_c:=\text{cor}^{-1}(c)\cap S$ the \textit{closed part of $S$}.
\end{defin}

Note that $\text{cor}$ being order preserving means that if $\text{cor}(S)=c$ then $\text{cor}(a)=c,\ \forall a\in S$.

\begin{defin}
    A \textit{2-function} between 2-sets $S$ and $S'$ is a based function $f:S_\ast\rightarrow S_\ast'$ such that $\text{cor}_S=\text{cor}_{S'}\circ f$. Note that there are no 2-functions between 2-sets of different colors, even if their underlying sets are the same.
\end{defin}

\begin{defin}
    Let $2\mathbb S_\text{inj}$ be the category whose objects are (possibly empty) finite 2-sets $S$ and whose morphisms are injective 2-functions. We will denote by $\text{sk}2\mathbb S_\text{inj}$ the skeleton of $2\mathbb S_\text{inj}$ whose object are $\underline{m,n}:=\{1_o,\dots,m_o,1_c,\dots,n_c\}_o$ and $\underline n:=\{1_c,\dots,n_c\}_c$ for $m,n\in \mathds N$. We also define $2\mathbb S_{o,\text{inj}}$ and $2\mathbb S_{c,\text{inj}}$ as the subcategories of open 2-sets and closed 2-sets respectively, and $\text{sk}2\mathbb S_{o,\text{inj}}$ and $\text{sk}2\mathbb S_{o,\text{inj}}$ its skeletons.
    
    We define the subcategories $2\mathbb S$, $\text{sk}2\mathbb S$, $2\mathbb S_{o}$, $2\mathbb S_{c}$, $\text{sk}2\mathbb S_{o}$ and $\text{sk}2\mathbb S_{c}$ as the subcategories containing only the bijective 2-functions.
\end{defin}

For every $S\in 2\mathbb S$ the 2-bijections $2\mathbb S_S:=2\mathbb S(S,S)$ have a group structure.

\begin{defin}
    Let $(\mathcal T,\otimes,\mathds 1)$ be a bicomplete symmetric monoidal category. A \textit{2$\mathbb S$-object in $\mathcal T$} is a functor $\mathcal O:\text{$2\mathbb S^\text{op}$}\rightarrow\mathcal T$ such that $\mathcal O$ is reduced, that is $\mathcal O(\underline 0)=\ast$ and $\mathcal O(\underline{0,0})=\ast$.
\end{defin}

For a $2\mathbb S$-object in $\mathcal T$ and $S\in 2\mathbb S$ we get an action of $2\mathbb S_{S}$ on $\mathcal O(S)$. A $2\mathbb S$-object in $\texttt{Top}$ is referred to as a $2\mathbb S$-space.

\begin{defin}
     A $2\mathbb S$-space $\mathcal O$ is $2\mathbb S$-free if the $2\mathbb S_{S}$ actions on $\mathcal O(S)$ are free.
\end{defin}

Again, the language of trees allows us to define a monad in the category of 2$\mathbb S$-object in $\mathcal T$.

\begin{defin}
    A 2-tree is a tree $T$ equipped with an order preserving function $\text{cor}_T:E_T\rightarrow 2$. An ordering on a 2-tree $T$ is a linear order of $V_T$.
    
    The coloring of a 2-tree $T$ induces a 2-set structure on $\text{in}(T)$ with $\text{cor}_{\text{in}(T)}(e)=\text{cor}_{T}(e)$ for $e\in\text{in}(T)$ and $\text{cor}_{\text{in}(T)}(\text{in}(T))=\text{cor}_T(e_0)$. A \textit{tree on a 2-set $S$} is a 2-tree $T$ equipped with a 2-bijection between $\text{in}(T)$ and $S$.
    
    We denote by $\mathbb T(2\mathbb S)$ the category whose objects are isomorphism classes of ordered 2-trees on 2-sets and whose morphisms are (not necessarily order preserving) isomorphisms of 2-trees. We also define $\mathbb T(\text{sk}2\mathbb S)$ the subcategory of trees on the sets of the form $\underline m=\{1,\dots,m\}$.
\end{defin}

For a 2-tree $T$ and vertex $v\in V_T$ the coloring of $T$ induces a 2-set structure on $\text{in}(v)$, with $\text{cor}_{\text{in}(v)}(e)=\text{cor}_{T}(e)$ for $e\in\text{in}(v)$ and $\text{cor}_{\text{in}(v)}(\text{in}(v))=\text{cor}_T(\text{out}(v))$.

Just like for trees, there is a grafting operation on 2-trees.

\begin{defin}
    For each $2\mathbb S$-object $\mathcal O$ define the $2\mathbb S$-object
    \begin{align*}
        \mathbb T_2\mathcal O: 2\mathbb{S}^{\text{op}}&\rightarrow \mathcal T\\
        S&\mapsto \text{Coeq}\left(\coprod_{T'\rightarrow T}\bigotimes_{v\in V_{T'}}\mathcal O(\text{in}(v))\rightrightarrows\coprod_{T\in \mathbb T(S)}\bigotimes_{v\in V_T}\mathcal O(\text{in}(v))\right)
    \end{align*}
    For $S=\underline{1,0},\underline 1$ the summand associated to the open and closed 2-tree without vertices is $\mathds 1$.
    
    As before this defines an endofunctor $\mathbb T_2$ of $\mathcal T^{2\mathbb S^{\text{op}}}$ that is the underlying functor of a monad.
\end{defin}

\begin{defin}
    A 2-operad $\mathcal O$ in $\mathcal T$ is a $\mathbb T_2$-algebra in $\mathcal T^{2\mathbb S^{\text{op}}}$.
    
    The 2-operad structure of $\mathcal O$ is completely defined by morphisms
$$
\eta_o^{\mathcal O}:\mathds 1\rightarrow \mathcal O(\underline{1,0}),\ \ \ \ \ \eta_c^{\mathcal O}:\mathds 1\rightarrow \mathcal O(\underline{1})
$$
and, for each $P$ linearly ordered 2-set and family of finite 2-sets $S_p$ indexed by $P$ with $\text{cor}(S_p)=\text{cor}(p)$ a morphism
$$
    \mu^{\mathcal O}_\nu:\mathcal O(P)\otimes \left(\bigotimes_{p\in P}\mathcal O(S_p)\right)\rightarrow\mathcal O(\cup_{p\in P}S_p)
$$
Again, these morphisms satisfy associative, unit and equivariance laws.
\end{defin}

Any 2-operad in $\mathcal T$ defines a monad in $\mathcal T_\ast^2$.

If $\mathcal O$ is a 2-operad we can as before extend the underlying functor on $2\mathbb S^{\text{op}}$ to a functor on $2\mathbb S_{\text{inj}}^{\text{op}}$, and any $(X_o,X_c)\in \mathcal T_\ast^2$ defines a functor:
\begin{align*}
    (X_o,X_c)^-:2\mathbb S_{\text{inj}}&\rightarrow \mathcal T\\
    S&\mapsto X_o^{S_o}\times X_c^{S_c}
\end{align*}

We can then define for any 2-operad $\mathcal O$ the monad
\begin{align*}
    O:\mathcal T^2_\ast&\rightarrow \mathcal T^2_\ast\\
    (X_o,X_c)&\mapsto \left(\int^{\text{sk}2\mathbb S_{o,\text{inj}}}\mathcal O(\underline{m,n})\times (X_o,X_c)^{\underline{m,n}},\int^{\text{sk}2\mathbb S_{c,\text{inj}}}\mathcal O(\underline n)\times X_c^{\underline{n}}\right)
\end{align*}
with the unit given by $(x_o,x_c)\mapsto ([\eta^{\mathcal O}_o,x_o],[\eta^{\mathcal O}_c,x_c])$ and multiplication induced by the 2-operad multiplications $\mu^{\mathcal O}$.

\begin{defin}
    For a 2-operad $\mathcal O$ an $\mathcal O$-algebra is an $O$-algebra for the associated monad $O$.
\end{defin}

For every 2-operad $\mathcal O$ and $(X_o,X_c)\in\texttt{Top}_\ast^2$ the space $O(X_o,X_c)$ admits a natural double filtration. Let $F^{k,l}\text{sk}2S_{\text{inj}}$ be the full subcategory of $\text{sk}2S_{\text{inj}}$ containing only the closed 2-sets $\underline n$ and the open 2-sets $\underline{m,n}$ with $m\leq k$ and $n\leq l$. We can then define $F^{k,l}\mathcal O(X_o,X_c)$ as the image of the natural inclusion of
$$
    \left(\int^{F^{k,l}\text{sk}2\mathbb S_{o,\text{inj}}}\mathcal O(\underline{m,n})\times (X_o,X_c)^{\underline{m,n}},\int^{F^l\text{sk}2\mathbb S_{c,\text{inj}}}\mathcal O(\underline n)\times X_c^{\underline{n}}\right)
$$
in $O(X_o,X_c)$. This of course induces a filtration $F^{k}\mathcal O(X_o,X_c):=\cup_{l\in\mathds N}F^{k,l}\mathcal O(X_o,X_c)$.

\subsection{Swiss-cheese 2-operads}

The swiss cheese 2-operads $\mathcal{SC}^N$ were first introduced by Voronov in \cite{Vo99}. For $N<\infty$ define $\mathcal{SC}^N(S)$ as the set of configurations of little $N$-cubes $d_S$ indexed by the elements of the 2-set $S$ such that the images of the $d_a$ for $a\in S$ are pairwise disjoint and for $a\in S_o$ the image of $\{0\}\times I^{N-1}$ by $d_a$ are contained in $\{0\}\times I^{N-1}$. Note that for closed 2-sets this is the same definition as that for the little $N$-cubes operad $\mathcal{C}^N$. We can then topologize $\mathcal{SC}^N(S)$ as a subspace of ${I^N}^{\coprod_{S}  I_a^{N}}$.

The $2\mathbb S$-space structural maps induced by the bijections of $S$ switch the indices of the $d_a$. The 2-operad structure is given by composition of the little cubes. The degeneracies delete the little $N$-cubes indexed by the elements that are not in the image of the inclusions $S\rightarrow S'$.

We have natural 2-operad inclusions
\begin{align*}
    \sigma^N_S:\mathcal{SC}^N(S)&\rightarrow\mathcal{SC}^{N+1}(S)\\
     \coprod_{a\in S} d_a&\mapsto\coprod_{a\in S} d_a \times1_I
\end{align*}
and we define $\mathcal{SC}^\infty:=\varinjlim \mathcal{SC}^N$, with the topology of the union.

\begin{prop}
    The images of the $N$-fold 2-loop space functors $\Omega_2^N$ are naturally $\mathcal{SC}^N$-spaces.
\end{prop}

\textbf{Proof:} First assume $N<\infty$. Let $\theta_2^N:SC^{N}\Omega^N_2\rightarrow \Omega^N_2$ be the natural transformation defined by
\begin{align*}\small
    \theta_{2}^N:SC^N\Omega^N_2\iota&\rightarrow \Omega^N_2\iota\\
    \left[d_{\underline{m,n}},((\alpha,\gamma)^{\underline m},\delta^{\underline n})\right]&\mapsto\\
    \noalign{\centering$\left(s\mapsto\begin{cases}
        \left(\alpha^i(d_i^{-1}(0,s)),s'\mapsto\begin{cases}\gamma^i(d_i^{-1}(0,s))(s'),&(s',s)\in d_i(I^N)\\
        \delta^j(d_j^{-1}(s',s)),&(s',s)\in d_j(I^N)\\
        \ast,&(s',s)\not\in d_{\underline{m,n}}(\coprod_{\underline{m,n}}I^N)
        \end{cases}\right),&(0,s)\in d_i(I^N)\\
        \left(\ast,s'\mapsto\begin{cases}\delta^j(d_j^{-1}(s',s)),&(s',s)\in d_j(I^N)\\
        \ast,&(s',s)\not\in d_{\underline{m,n}}(\coprod_{\underline{m,n}}I^N)
        \end{cases}\right),&(0,s)\not\in d_{\underline m}(\coprod_{\underline{m}}I^N)
    \end{cases}\right)$}
    \left[d_{\underline{n}},\delta^{\underline n}\right]&\mapsto\left(t\mapsto\begin{cases}
        \delta^j(d_j^{-1}(t)),&t\in d_j(I^N)\\
        \ast, &t\not\in d_{\underline n}(\coprod_{\underline{n}}I^N)
    \end{cases}\right)
\end{align*}

It is straightforward to verify that $\theta_{2}^N$ is compatible with the $2$-operad structural morphisms of $\mathcal{SC}^N$.

The action of $\mathcal{SC}^\infty$ on the image of $\Omega_2^\infty$ is induced from the finite cases. By definition an element $d_{\underline{m,n}}\in \mathcal{SC}^\infty(\underline{m,n})$ belongs to some $\mathcal{SC}^N(\underline{m,n})$ with $N<\infty$. We can then use the above construction for the finite case to define the structural natural transformation $\theta^\infty$.$\blacksquare$

\subsection{Model structure on 2-operads and algebras over 2-operads}

In order to describe the model structure on colored operads and algebras Berger and Moerdijk use the notion of coalgebra interval objects in \cite{BM07}.

\begin{defin}
    In a closed monoidal model category $\mathcal T$ an \textit{interval object} is an object $I\in\mathcal T$ equipped with a factorization  into a cofibration followed by a weak equivalence
    $$
        \mathds 1 \coprod \mathds 1 \xrightarrow{0\coprod 1} I \xrightarrow{\varepsilon} \mathds 1
    $$
    of the codiagonal morphism on the unit $\mathds 1$, and with a monoid structure $(I,\vee,0)$ such that $1$ is an absorbing element, that is $x\vee 1=1=1\vee x$, and $\varepsilon$ is a counit, that is $\varepsilon(x\vee y)=\varepsilon(x)\varepsilon(y)$ and $\varepsilon\circ 0=Id_{\mathds 1}=\varepsilon \circ 1$.
    
    A \textit{(cocommutative) coalgebra interval object} is an interval object $I$ equipped with a (cocommutative) comonoid structure such that the morphisms in the factorization of the codiagonal are comonoid morphisms.
\end{defin}

In $\texttt{Top}$ the interval $I=[0,1]$ with the maximum operation as the monoid structure and the codiagonal comonoid structure is a cocommutative coalgebra interval object.

From \cite{BM07} we get the following conditions for the existence of a model structure on 2-operads and algebras over 2-operads.
 
\begin{prop}
    Let $\mathcal T$ be a cofibrantly generated closed symmetric model category with a symmetric monoidal fibrant replacement functor and a coalgebra interval object. Then there is a cofibrantly generated model structure on the category of 2-operads in $\mathcal T$, in which a morphism $f:\mathcal O\rightarrow \mathcal P$ is a weak equivalence (resp. fibration) if and only if for each $S\in 2\mathbb S$ the morphism $f_S:\mathcal O(S)\rightarrow \mathcal P(S)$ is a weak equivalence (resp. fibration) in $\mathcal T$.
    
    Moreover if the coalgebra interval object is cocommutative then for every 2-operad $\mathcal O$ the $(U\dashv O)$-adjunction transfers the model structure from $\mathcal T_\ast$ to $\mathcal O[\mathcal T]$.
\end{prop}

If $(X_o,X_c)\in\texttt{Top}_\ast^2$ is a cofibrant object weakly equivalent to $\Omega^N_2\iota$ for some $\iota:B\rightarrow Y\in\texttt{Top}_\ast^\shortrightarrow$ it is not generally the case that $(X_o,X_c)$ is a $\mathcal{SC}^N$-space. This is a consequence of the fact that $\mathcal{SC}^N$ is not a cofibrant 2-operad. The theorem \cite[Thm. 3.5]{BM03} implies that $(X_o,X_c)$ does have a $\overline{\mathcal{SC}^N}$-space structure for any cofibrant resolution $\overline{\mathcal{SC}^N}$ of $\mathcal{SC}^N$.

\begin{prop}\label{Act M.3.4}
    Let $\psi:\mathcal O\rightarrow\mathcal O'$ be a weak equivalence of $2\mathbb S$-free 2-operads and $(X_o,X_c)\in\texttt{Top}^2_\ast$, then $\psi_{(X_o,X_c)}:O(X_o,X_c)\rightarrow O'(X_o,X_c)$ induces an isomorphism on integral homology.
\end{prop}

\textbf{Proof:} The proof is similar to the one for \cite[Prop. 3.4]{Ma72}.$\blacksquare$

The following will allow us to define the left derived functor of the bar construction on certain cases.

\begin{lema}\label{LeftDerB}
    If $\mathcal T$ is a closed monoidal model category, $\mathcal O$ is a cofibrant 2-operad and $S$ is an $O$-functor in $\mathcal A$ that is a weak Quillen left adjoint, then
    $$
        B_\bullet(S,O,-):\mathcal O[\mathcal T]\rightarrow \mathcal A^{\Delta^{\text{op}}}
    $$
    admits a left derived functor.
\end{lema}

\textbf{Proof:} The coend over the tensor is a left Quillen bifunctor \cite[A.2.9.26]{Lu09}, so with $\mathcal O$ cofibrant the monad $O$ is left Quillen adjoint, which means it presererves (trivial) cofibrations. Also the unit $\eta$ is a cofibration on cofibrant objects. Therefore $B_\bullet(S,O,X)$ is cofibrant if $X$ is cofibrant. If $f:X\rightarrow Y$ is a weak equivalence between cofibrant objects then $SO^qf$ is a weak equivalence for all $q$, and therefore $B(1,1,f)$ is a weak equivalence.$\blacksquare$

Since geometric realization is a left Quillen adjoint functor when $\mathcal A$ is a topological category the above lemma also implies that $B(S,O,-)$ admits a left derived functor.

\section{Relative recognition principle}

In 4.1 the approximation theorem is reviewed, and we see that there is a natural transformation $\alpha^N_2$ between $SC^N$ and $\Omega^N_2\Sigma^N_2$ that is a weak equivalence on certain $(X_o,X_c)$. In 4.2 some compatibility results between the geometric realization, $\Omega^N_2$, $\Sigma^N_2$ and $\mathcal{SC}^N$ are stated without proof (the arguments are analogous to results in \cite[12]{Ma72}). In 4.3 and 4.4 the relative recognition principle is proved in the $N<\infty$ and $N=\infty$ cases respectively.

\subsection{Approximation theorem}

Let $\alpha^N$ be the composition of natural transformations
$$
    C^N\xrightarrow{C^N \eta^N}C^N\Omega^N\Sigma^N\xrightarrow{\theta^N_{\Sigma^N}}\Omega^N\Sigma^N
$$
which is explicitly given by
    \begin{align*}
    \alpha^N_X:C^NX&\rightarrow \Omega^N \Sigma^N X\\
    \left[d_{\underline{n}},x^{\underline{n}}\right]&\mapsto \left(t\mapsto\begin{cases}
        [x^i,d_i^{-1}(t)],&t\in d_i(I^N)\\
        \ast,&t\not\in d_{\underline n}\left(\coprod_{\underline n}I^N\right)
    \end{cases}\right)
\end{align*}

We need some definitions in order to state the approximation theorem. See \cite{Ma74} for details.

\begin{defin}
    An $\mathcal H$-space is a space $X\in\texttt{Top}$ equipped with a map $\mu:X\times X\rightarrow X$ and an element $e\in X$ such that both $\mu(e,-)$ and $\mu(-,e)$ are homotopy equivalent to the identity.
\end{defin}

For an $\mathcal H$-space $(X,\mu,e)$ the homology groups $H_\bullet(X;k)$ for any commutative coefficients ring $k$ equipped with the  Pontryagin product $\mu_\ast$ and the unit $[e]$ is a graded $k$-algebra.

\begin{defin}
    An $\mathcal H$-space $X$ is homotopy associative if $\mu(-,\mu(-,-))$ is homotopy equivalent to $\mu(\mu(-,-),-)$.
\end{defin}

The $k$-algebra structure on $H_\bullet(X;k)$ for a homotopy associative $\mathcal H$-space $X$ is associative.

For every $d_{\underline 2}\in\mathcal C^N(\underline 2)$ and $\mathcal C^N$-space $(X,\xi)$ the map $\xi([d_{\underline 2},(-,-)])$ and the basepoint $\ast_X$ endows $X$ with a homotopy associative $\mathcal H$-space structure.

\begin{defin}
    An $\mathcal H$-space $X$ is admissible if it is homotopy associative and $\mu(x,-):$ and $\mu(-,x)$ are homotopy equivalent $\forall\ x\in X$.
\end{defin}

The $k$-algebra structure on $H_\bullet(X;k)$ for an admissible $\mathcal H$-space $X$ is associative and commutative.

For $1<N\leq \infty$ the $\mathcal H$-space structures on the $\mathcal C^N$-spaces are admissible.

\begin{defin}
    An $\mathcal H$-space $X$ is \textit{grouplike} if $(\pi_0 X,\mu_\ast, [e])$ has a group structure.
\end{defin}

The $\mathcal H$-space structures on $N$-fold loop spaces induced by the $\mathcal C^N$-structure are grouplike.

For $(X_o,X_c)\in\mathcal{SC}^N[\texttt{Top}]$ we have that $X_c$ is a $\mathcal C^N$-space and $X_o$ is a $\mathcal C^{N-1}$-space. For $N>2$ we say that $(X_o,X_c)$ is grouplike if the $\mathcal H$-space structures on both $X_o$ and $X_c$ are grouplike. We denote by $\mathcal{SC}^N[\texttt{Top}]_{\text{grp}}$ the subcategory of grouplike $\mathcal{SC}^N$-spaces.

\begin{defin}
    A homological group completion is an $\mathcal H$-map $g:X\rightarrow G$ between admissible $\mathcal H$-spaces such that $G$ is grouplike and for every commutative ring $k$
    $$
        \bar g_\ast:H_\bullet(X,k)[\pi_0 X^{-1}]\rightarrow H_\bullet (G,k)
    $$
    is an isomorphism.
\end{defin}

For every grouplike homotopy associative $\mathcal H$-space $X$ we have a homotopy equivalence between $X$ and $X_e\times\pi_0X$, where $X_e$ is the connected component containing $e$ \cite[Lemma I.4.6]{CLM76}. Since every $\mathcal H$-space $X$ is simple, that is $\pi_1 X$ is abelian and acts trivially on $\pi_q X$ for all $q$, the dual Whitehead theorem for connected $\mathcal H$-spaces implies that a group completion of a grouplike admissible $\mathcal H$-space is a weak equivalence.

\begin{teor}[Approximation theorem] If $X\in\texttt{Top}_\ast$ is connected,
    $\alpha^N_X$ is a weak equivalence. If $1<N\leq \infty$ then $\alpha^N_X$ is a homological group completion for all $X\in\texttt{Top}_\ast$.
\end{teor}

\begin{cor}\label{CorApprox}
    If $X\in\mathcal{C}^N[\texttt{Top}]$ is connected, then $\alpha^N_X$ is a weak equivalence. If $1<N\leq \infty$ and $X\in\mathcal{C}^N[\texttt{Top}]$ is grouplike then $\alpha^N_X$ is a weak equivalence.
\end{cor}

We will prove a relative version of this corollary. Define $\alpha_2^N$ as the composition of natural transformations
$$
    SC^N\xrightarrow{SC^N\eta_2^N}SC^N\Omega_2^N\Sigma_2^N\xrightarrow{\theta^N_{ \Sigma_2^N}}\Omega_2^N\Sigma_2^N.
$$
which is explicitly given by
\begin{align*}\small
    \alpha^N_{2(X_o,X_c)}:SC^N (X_o,X_c)&\rightarrow \Omega_2^N\Sigma_2^N(X_o,X_c)\\
    \left[d_{\underline{m,n}}, x^{\underline{m,n}}\right]&\mapsto\\
    \noalign{\centering$\left(s\mapsto\begin{cases}
        \left([x^{i},d_{i}^{-1}(0,s)],s'\mapsto\begin{cases}[x^{a},d_{a}^{-1}(s',s)],&(s',s)\in d_a(I^N)\\
        \ast,&(s',s)\not\in d_{\underline{m,n}}(\coprod_{\underline{m,n}}I^N)
        \end{cases}\right),&(0,s)\in d_i(I^N)\\
        \left(\ast,s'\mapsto\begin{cases}[x^{j},d_{j}^{-1}(s',s)],&(s',s)\in d_j(I^N)\\
        \ast,&(s',s)\not\in d_{\underline{m,n}}(\coprod_{\underline{m,n}}I^N)
        \end{cases}\right),&(0,s)\not\in d_{\underline m}(\coprod_{\underline{m}}I^N)
    \end{cases}\right)$}
    \left[d_{\underline{n}},x^{\underline n}\right]&\mapsto\left(t\mapsto\begin{cases}
        [x^j,d_j^{-1}(t)],&t\in d_j(I^N)\\
        \ast, &t\not\in d_{\underline n}(\coprod_{\underline{n}}I^N)
    \end{cases}\right)
\end{align*}

Note that $\alpha^N_c$ is the composition of $\alpha^N$ with the $\Omega^N$ image of the inclusion $\Sigma^NX_c\hookrightarrow \Sigma^{N-1}(X_o\wedge I)\vee \Sigma^NX_c$, which is the inclusion of a deformation retract and therefore a weak equivalence.

Let $\pi_o, \pi_c:\texttt{Top}_\ast^2\rightarrow \texttt{Top}_\ast$ be the projections on the open and closed coordinate respectively.

The functors $\pi_c SC^N$ and $\pi_o SC^N$ are $SC^N$-functors, since we can define $\lambda^{\pi_o SC^N}:=\pi_o\mu^{SC^N}$ and $\lambda^{\pi_c SC^N}:=\pi_c\mu^{SC^N}$. Note that $\pi_c SC^N=C^N\pi_c$ and $\pi_c\mu^{SC^N}=\mu^{C^N}_{\pi_c}$.

The functor $C^{N-1}\pi_o$ has the structure of a $SC^N$-functor. In order to build the $SC^N$-functor structural map $\lambda^{C^{N-1}\pi_o}:C^{N-1}\pi_o SC^N\rightarrow C^{N-1}$ first note that for any $\mathbb S$-object $\mathcal O$ we can extend it to a $2\mathbb S$-object by
\begin{align*}
    \mathcal O:2\mathbb S&\rightarrow\mathcal T\\
    S&\mapsto \mathcal O(S_o)
\end{align*}
    We can then define the natural transformation $p_o:\mathcal{SC}^N\rightarrow \mathcal C^{N-1}$ as
\begin{align*}
    p_o:\mathcal{SC}^N(S)&\rightarrow C^{N-1}(S_o)\\
    d_S&\mapsto
    d_{S_o}\restriction_{\{0\}\times I^{N-1}}
\end{align*}
which we can use to define the $SC^N$-functor structural map as follows,
\begin{align*}
    \lambda^{C^{N-1}\pi_o}:C^{N-1}\pi_o SC^N(X_o,X_c)&\rightarrow C^{N-1}X_o\\
    \left[d_{\underline{m'}},\prod_{i\in\underline{m'}}[d_{\underline{m_i,n_i}},x^{\underline{m_i,n_i}}]\right]&\mapsto \left[\mu^{\mathcal C^{N-1}}\left(d_{\underline{m'}}\otimes \left(\bigotimes_{i\in {\underline{m'}}}p_od_{\underline{m_i,n_i}} \right)\right),x^{\coprod_{i\in\underline{m'}}\underline{m_i}}\right]
\end{align*}

The natural transformation $p_o$ also allows us to define the following map between $SC^N$-functors:
\begin{align*}
    p^N_o: \pi_o SC^N(X_o,X_c)&\rightarrow C^{N-1}X_o\\
    \left[d_{\underline{m,n}},x^{\underline{m,n}}\right]&\mapsto\left[p_od_{\underline{m,n}},x^{\underline m}\right]
\end{align*}

Note that $\pi_o\mathcal{SC}^N(X_o,X_c)$ is a $\mathcal C^{N-1}$-space and that $p_o^N$ is also a $\mathcal C^{N-1}$-map, and therefore when $N>1$ it is an $\mathcal H$-map and when $N>2$ it is an $\mathcal H$-map between admissible $\mathcal H$-spaces. We now prove $p_o^N$ is a quasifibration.

\begin{defin}
    A map $p: E \rightarrow B$ in $\texttt{Top}$ is a \textit{quasifibration} if the natural inclusions 
    \begin{align*}
        p^{-1}(b)&\rightarrow \{(e,\gamma)\in E\times_BB^I\ |\ \gamma(0)=b,\gamma(1)=e\}\\
        f&\mapsto (f,t\mapsto p(f))
    \end{align*}
        
    are weak equivalences for all $b\in B$
    
    We say that a subspace $U\subset B$ is distinguished if $p:p^{-1}(U)\rightarrow U$ is a quasifibration.
\end{defin}

From \cite[2.7]{Ma90} we get the following criterion for a map to be a quasifibration.

\begin{prop}\label{CondQsiFib2}
     Let $p: E \rightarrow B$ be a map of filtered spaces such that $F^kE = p^{-1}F^kB$ for each $k \geq 0$. If for each $k\geq 1$ the map $p: F^kE \rightarrow F^kB$ is obtained by passage to pushouts from a commutative diagram of the form
     $$
        \xymatrix{
            F^{k-1}E\ar[d]^p&D_k\ar[l]_{g_k}\ar[d]^{q_k}\ar[r]^{j_k}&E_k\ar[d]^{p_k}\\
            F^{k-1}B&\ar[l]^{f_k}A_k\ar[r]_{i_k}&B_k
        }
     $$
     such that the following conditions hold:
     \begin{enumerate}[i)]
         \item $F^0B$ is distinguished;
         
         \item $p_k:E_k\rightarrow B_k$ is a fibration;
         
         \item $i_k:A_k\rightarrow B_k$ and $j_k:D_k\rightarrow E_k$ are inclusions of NDR-pairs;
         
         \item The right square is a pullback;
         
         \item $g_k:q_k^{-1}(a) \rightarrow p^{-1}(f_k(a))$ is a weak equivalence for all $a \in A_k$;
     \end{enumerate}
     then each $F^kB$ is distinguished and $p: E \rightarrow B$ is a quasifibration.
\end{prop}

\begin{teor}
If $X_o$ is well-pointed, the map $p_o^N$ is a quasifibration with fiber $C^N(X_c)$.
\end{teor}

\textbf{Proof:} We'll show that $p_o^N$ and the natural filtrations on $\mathcal C^{N-1}X_o$ and $\mathcal{SC}^N(X_o,X_c)$ from section 3 satisfy the conditions of proposition \ref{CondQsiFib2}.

Let $\mathbb S_{\underline k}$ be the full subcategory of $\text{sk}\mathbb S$ containing only $\underline k$ and $2\mathbb S_{\underline k}$ be the full subcategory of $\text{sk}2\mathbb S_{\text{inj}}$ containing only the open 2-sets of the form $\underline{k,n}$.

We can then define
$$
    C^{N-1}_{k}X_o:=\int^{\mathbb S_{\underline k}}\mathcal C^{N-1}(\underline{k})\times X_o^{\underline{k}}
$$
and
$$
    SC^{N}_{k}(X_o,X_c):=\int^{2\mathbb S_{\underline k}} \mathcal{SC}^{N}(\underline{k,n})\times (X_o,X_c)^{\underline{k,n}}
$$

Define also
$$
    A^N_k:=\{[d_{\underline k},x^{\underline k}]\in C^{N-1}_{k}X_o\ |\ \exists i\in\underline k:x^i=\ast_o\}
$$
and
$$
    D^N_k:=\{[d_{\underline{k,n}},x^{\underline{k,n}}]\in SC^{N}_{k}(X_o,X_c)\ |\ \exists i\in\underline k:x^i=\ast_o\}
$$

The maps $p^N_o:F^k\pi_o\circ SC(X_o,X_c)\rightarrow F^kC^{N-1}X_o$ are then the pushouts of

$$
\xymatrix{
            F^{k-1}\pi_o SC^N(X_o,X_c)\ar[d]^{p^N_o}&D^N_k\ar[l]_(0.3){g_k}\ar[d]^{q_k}\ar[r]^(0.4){j_k}&SC^{N}_{k}(X_o,X_c)\ar[d]^{p_k}\\
            F^{k-1}C^{N-1}X_o&\ar[l]^(0.4){f_k}A^N_k\ar[r]_(0.4){i_k}&C^{N-1}_kX_o
        }
$$
where $i_k$ and $j_k$ are the inclusions of subspaces, $g_k$ and $f_k$ are induced by the degeneracy of the little cubes with the same index as the $x^i$ that are equal to the base point and $q_k$ and $p_k$ are defined in the same way as $p^N_o$. These diagrams satisfy the conditions in proposition \ref{CondQsiFib2}:

\begin{enumerate}[i)]
\item $F^{0}C^{N-1}X_o=\{\ast\}$. Since every space is fibrant and every fibration is a quasifibration $F^{0}C^{N-1}X_o$ is distinguished;

\item By definition $p_k$ is a fibration if for every commutative square
$$
    \xymatrix{
        I^q\ar@{^{(}->}[d]_{\iota_0^q}\ar[r]^(0.35){[d_{\underline{k,l}},x^{\underline{k,l}}]}&SC^N_k(X_o,X_c)\ar[d]^{p_k}\\
        I^q\times I\ar@{.>}[ur]|{\tilde H}\ar[r]_(0.45){[\delta_{\underline{k}},\xi^{\underline{k}}]}&C^{N-1}_kX_o
    },
$$
there is a lift $\tilde H$ that makes the diagram commute.

Remember from the beginning of section 3.4 that a little cube $d$ in $\mathcal{SC^{N}}(\underline{m,n})$ or $\mathcal{C^{N}}(\underline{m})$ is of the form $d(s)=M_ds+C_d$ for some diagonal $N\times N$ matrix $M_d$ with strictly positive entries and some $C_d\in I^N$. Therefore for any $(d_{\underline{m,n}},\delta_{\underline m})\in \mathcal{SC}^{N}(\underline{m,n})\times_{\mathcal{C}^{N-1}(\underline{m})}{\mathcal{C}^{N-1}(\underline{m})}^{I}$ and $v\in(0,1]$ we can define
\begin{align*}
    \gamma^{\nu}_{d_{\underline{m,n}},\delta_{\underline{m}}}:I&\rightarrow {I^N}^{\coprod_{\underline{m,n}}I^N_a}\\
    t&\mapsto \left(s_a\mapsto\begin{cases}\begin{cases}
    \begin{bmatrix}(\frac{2v-t}{2v})m^1_{d_{a}}&0\\
    0&M_{\delta_{a}(t)}\end{bmatrix}s_a+\begin{bmatrix}(\frac{2v-t}{2v})c^1_{d_{a}}\\
    C_{\delta_{a}(t)}\end{bmatrix},&0\leq t\leq v\\
    \begin{bmatrix}\frac{1}{2}m^1_{d_{a}}&0\\
    0&M_{\delta_{a}(t)}\end{bmatrix}s_a+\begin{bmatrix}\frac{1}{2}c^1_{d_{a}}\\
    C_{\delta_a(t)}\end{bmatrix},&v\leq t\leq1
    \end{cases},&\text{cor}(a)=o;\\
    \begin{cases}\frac{2v-t}{2v}M_{d_{a}}s_a+\begin{bmatrix}\frac{(2v-t)c^1_{d_{a}}+t}{2v}\\
    \pi_{I^{N-1}} C_{d_{a}}\end{bmatrix},&0\leq t\leq v\\
    \frac{1}{2}M_{d_{a}}s_a+\begin{bmatrix}\frac{1+c^1_{d_{a}}}{2}\\
    \pi_{I^{N-1}} C_{d_{a}}\end{bmatrix},&v\leq t\leq1
\end{cases},&\text{cor}(a)=c.
    \end{cases}\right),
\end{align*}
i.e. $\gamma^{\nu}_{d_{\underline{m,n}},\delta_{\underline{m}}}$ is the path in ${I^N}^{\coprod_{\underline{m,n}}I^N_a}$ such that $\gamma^{\nu}_{d_{\underline{m,n}},\delta_{\underline{m}}}(0)=d_{\underline{m,n}}$, $p_o\gamma^{\nu}_{d_{\underline{m,n}},\delta_{\underline{m}}}=\delta_{\underline{n}}$, the heights in the first coordinate of the open little cubes are linearly halved in the interval $[0,\nu]$ and remain constant in the interval $[\nu,1]$, and such that the side lengths in all coordinates of the closed little cubes and the distance in the first coordinate from the center of the little cube to 1 are also linearly halved in the interval $[0,\nu]$ and remain constant in the interval $[\nu,1]$.

\begin{figure}[h!]
\centering\begin{tikzpicture}

\draw (0,0)--(0,2)--(2,2)--(2,0)--(1.7,0)--(1.7,1.5)--(1.3,1.5)--(1.3,0)--(0.9,0)--(0.9,0.8)--(0.5,0.8)--(0.5,0)--(0,0);

\draw (0.2,0.9) rectangle (0.7,1.2);

\draw (0.8,1.2) rectangle (1.2,1.8);

\draw (0,0)--(1,-3)--(1.1,-3);

\draw (2,0)--(3,-3)--(2.9,-3);

\draw (1.4,-3)--(2,-3);

\draw (0.5,0) .. controls (0.2,-1) and (1.3,-2.5) .. (1.1,-3);

\draw (0.9,0) .. controls (1.2,-0.7) and (1.6,-2) .. (1.4,-3);

\draw (1.3,0) .. controls (0.9,-0.5) and (2.4,-1.5) .. (2,-3);

\draw (1.7,0) .. controls (2,-0.5) and (2.5,-2) .. (2.9,-3);

\draw[dashed] (0.666,-2)--(2.666,-2);

\draw (0.45,-2) node {$\nu$};

\end{tikzpicture}\begin{tikzpicture}

\draw (0.5,0)--(0,0)--(1,-3)--(1.1,-3);

\draw (1.7,0)--(2,0)--(3,-3)--(2.9,-3);

\draw (1.4,-3)--(2,-3);

\draw (0.9,0)--(1.3,0);

\draw (0.5,0) .. controls (0.2,-1) and (1.3,-2.5) .. (1.1,-3);

\draw (0.9,0) .. controls (1.2,-0.7) and (1.6,-2) .. (1.4,-3);

\draw (1.3,0) .. controls (0.9,-0.5) and (2.4,-1.5) .. (2,-3);

\draw (1.7,0) .. controls (2,-0.5) and (2.5,-2) .. (2.9,-3);

\draw[fill=white] (0.333,-1)--(0.333,1)--(2.333,1)--(2.333,-1)--(2.165,-1)--(2.165,0.125)--(1.59,0.125)--(1.59,-1)--(1.24,-1)--(1.24,-0.4)--(0.533,-0.4)--(0.533,-1)--(0.333,-1)

 (0.5955,0.175)-- (0.9705,0.175)-- (0.9705,0.4) -- (0.5955,0.4)--(0.5955,0.175)

 (1.182,0.4) -- (1.482,0.4)-- (1.482,0.85) -- (1.182,0.85)-- (1.182,0.4);

\draw[dashed] (0.666,-2)--(2.666,-2);

\draw (0.45,-2) node {$\nu$};

\end{tikzpicture}\begin{tikzpicture}

\draw (0.5,0)--(0,0)--(1,-3)--(1.1,-3);

\draw (1.7,0)--(2,0)--(3,-3)--(2.9,-3);

\draw (1.4,-3)--(2,-3);

\draw (0.9,0)--(1.3,0);

\draw (0.5,0) .. controls (0.2,-1) and (1.3,-2.5) .. (1.1,-3);

\draw (0.9,0) .. controls (1.2,-0.7) and (1.6,-2) .. (1.4,-3);

\draw (1.3,0) .. controls (0.9,-0.5) and (2.4,-1.5) .. (2,-3);

\draw (1.7,0) .. controls (2,-0.5) and (2.5,-2) .. (2.9,-3);

\draw[dashed] (0.666,-2)--(2.666,-2);

\draw (0.45,-2) node {$\nu$};

\draw[fill=white] (0.666,-2)--(0.666,0)--(2.666,0)--(2.666,-2)--(2.51,-2)--(2.51,-1.25)--(2,-1.25)--(2,-2)--(1.43,-2)--(1.43,-1.6)--(0.88,-1.6)--(0.88,-2)--(0.666,-2)

 (0.991,-0.55)-- (1.241,-0.55)-- (1.241,-0.4) -- (0.991,-0.4)--(0.991,-0.55)

 (1.566,-0.4) -- (1.766,-0.4)-- (1.766,-0.1) -- (1.566,-0.1)-- (1.566,-0.4);

\end{tikzpicture}\begin{tikzpicture}

\draw (0.5,0)--(0,0)--(1,-3)--(1.1,-3);

\draw (1.7,0)--(2,0)--(3,-3)--(2.9,-3);

\draw (1.4,-3)--(2,-3);

\draw (0.9,0)--(1.3,0);

\draw (0.5,0) .. controls (0.2,-1) and (1.3,-2.5) .. (1.1,-3);

\draw (0.9,0) .. controls (1.2,-0.7) and (1.6,-2) .. (1.4,-3);

\draw (1.3,0) .. controls (0.9,-0.5) and (2.4,-1.5) .. (2,-3);

\draw (1.7,0) .. controls (2,-0.5) and (2.5,-2) .. (2.9,-3);

\draw[dashed] (0.666,-2)--(2.666,-2);

\draw (0.45,-2) node {$\nu$};

\draw[fill=white] (1,-3)--(1,-1)--(3,-1)--(3,-3)--(2.9,-3)--(2.9,-2.25)--(2,-2.25)--(2,-3)--(1.4,-3)--(1.4,-2.6)--(1.1,-2.6)--(1.1,-3)--(1,-3)

 (1.324,-1.55)-- (1.574,-1.55)-- (1.574,-1.4) -- (1.324,-1.4)--(1.324,-1.55)

 (1.9,-1.4) -- (2.1,-1.4)-- (2.1,-1.1) -- (1.9,-1.1)-- (1.9,-1.4);

\end{tikzpicture}
\caption{Cross sections of $\gamma^{\nu}_{d_{\underline{m,n}},\delta_{\underline{m}}}$.}
\end{figure}

For some $(d_{\underline{m,n}},\delta_{\underline m})$ and $\nu$ it might be the case that $\gamma^{\nu}_{d_{\underline{m,n}},\delta_{\underline m}}\not\in \mathcal{SC}(\underline{m,n})^{I}$, because there might be $a,a'\in \underline{m,n}$ and $t\in I$ such that $a\neq a'$ and $\gamma^{\nu}_{d_{\underline{m,n}},\delta_{\underline m},a}(t)(\mathring I^N)\cap \gamma^{\nu}_{d_{\underline{m,n}},\delta_{\underline m},a'}(t)(\mathring I^N)\neq\emptyset$, but there is always some $\nu'\in(0,1]$ such that for $\nu\in(0,\nu']$ these intersections are empty. We can then define for each $\underline{m,n}\in\text{sk}2\mathbb S$ the map
\begin{align*}
    \nu_{m,n}:SC^{N}(\underline{m,n})\times_{C^{N-1}(\underline{m})}{C^{N-1}(\underline{m})}^{I}&\rightarrow (0,1]\\
    (d_{\underline{m,n}},\delta_{\underline m})&\mapsto \max\left\{v\in(0,1]\ |\ \gamma^{\nu}_{d_{\underline{m,n}},\delta_{\underline m}}\in \mathcal{SC}^{N}(\underline{m,n})^{I}\right\}
\end{align*}

We now prove that for any given commutative square as above we can always choose a $\bar \nu\in(0,1]$ such that $\gamma^{\bar\nu}_{d_{\underline{k,l}}(s),\delta_{\underline{k}}(s)}\in\mathcal{SC}^{N}(\underline{m,n})^{I}$ for all $s\in I^q$.

Let $F^l2\mathbb S_{k}$ be the full subcategory of $2\mathbb S_{k}$ containing only the open 2-sets of the form $\underline{k,n}$ with $n\leq l$. The natural inclusions of $\int^{F^l2\mathbb S_{k}} \mathcal{SC}^N(\underline{k,n})\times (X_o,X_c)^{\underline{k,n}}$ in $SC^N_k(X_o,X_c)$ gives us a filtration $F^lSC^N_k(X_o,X_c)$.

Any compact subspace of $SC^N_k(X_o,X_c)$ is contained in $F^lSC^N_k(X_o,X_c)$ for some $l\in\mathds N$. To see this, assume on the contrary that there is an infinite sequence of points $z_i \in K$ all lying in distinct $F^lSC^N_k(X_o,X_c)$. Consider the subset $S$ of all such points in $K$. In order to show that $S$ is closed, assume that $S\cap F^{l-1}SC^N_k(X_o,X_c)$ is closed. Then $S \cap F^lSC^N_k(X_o,X_c)$ contains at most one more point. The space $\mathcal{SC}_k^N(X_o,X_c)$ is weakly Hausdorff, so points are closed and therefore $S\cap F^{l}SC^N_k(X_o,X_c)$ is closed. It follows that $S$ is closed. The same argument shows that any subset of $S$ is closed, so $S$ has the discrete topology. Being a closed subset of a compact set, $S$ must be compact. Therefore, $S$ has to be finite, a contradiction.\footnote{The author thanks Eduardo Hoefel for providing the argument in this paragraph through private communication.}

Now the filtration $F^lSC^N_k(X_o,X_c)$ and the map $[d_{\underline{k,l}},x^{\underline{k,l}}]$ from the commutative diagram induce a filtration $F^lI^q:=[d_{\underline{k,l}},x^{\underline{k,l}}]^{-1}(F^lSC^N_k(X_o,X_c))$. We can then define
\begin{align*}
    \nu^l:\overline{F^lI^q-F^{l-1}I^q}&\rightarrow (0,1]\\
    s&\mapsto \nu_{k,l}(d_{\underline{k,l}}(s),\delta_{\underline{k}}(s))
\end{align*}
Since $I^q$ is compact $\overline{F^lI^q-F^{l-1}I^q}$ is compact, and therefore the image of $\nu^l$ has a positive minimum.

By the two previous observations we can define
$$
    \bar\nu:=\min\{\nu^l(s)\ |\ l\in\mathds N; s\in \overline{F^lI^q-F^{l-1}I^q}\}.
$$

This gives us a map
\begin{align*}
    \tilde H: I^q\times I&\rightarrow SC^N_k(X_o,X_c)\\
    (s,t)&\mapsto \left[\gamma^{\bar\nu}_{d_{\underline{k,l}}(s),\delta_{\underline{k}}(s)}(t), \left(\xi^{\underline k}(s,t),x^{\underline l}(s)\right)\right]
\end{align*}
which makes the diagram commute.

\item Since $X_o$ is well-pointed there are maps $u:X_o\rightarrow I$ and $H:X_o\times I\rightarrow X_o$ making $(\ast_o,X_o)$ an NDR-pair, and therefore we can define the maps
\begin{align*}
    u':C^{N-1}_kX_o&\rightarrow I\\
    [d_{\underline{k}},x^{\underline{k}}]&\mapsto \min\{u(x^i)\ |\ i\in \underline k\}
\end{align*}
and
\begin{align*}
    H':C^{N-1}_kX_o\times I&\rightarrow C^{N-1}_kX_o\\
    ([d_{\underline{k}},x^{\underline{k}}],t)&\mapsto \left[d_{\underline{k}},\prod_{i\in\underline{k}}H(x^{i},t)\right]
\end{align*}
which gives us that $(A^N_{k},C^{N-1}_kX_o)$ is an NDR-pair. That $(D^N_{k},SC^{N}_k(X_o,X_c))$ is an NDR-pair follows by an analogous argument.

\item It is trivial to check that the right square is a pullback.

\item Fix $[d_{\underline k},x^{\underline k}
]\in A^N_k$ and define the subspaces
$$
    P^N_k:=\left\{[d_{\underline{k,n}},x^{\underline{k,n}}]\in q_k^{-1}([d_{\underline k},x^{\underline k}])\ \left|\ \begin{matrix*}[l]\forall i\in\underline k,\  d_i(I)\subset\left[0,\frac{1}{2}\right]\times I^{N-1}; \\\forall j\in\underline n,\ d_j(I)\subset\left[\frac{1}{2},1\right]\times I^{N-1}.\end{matrix*}\right.\right\}
$$
and
$$
    Q^N_k:=\left\{[d_{\underline{k-1,n}},x^{\underline{k-1,n}}]\in p_o^{-1}(f_k([d_{\underline k},x^{\underline k}]))\ \left|\ \begin{matrix*}[l]\forall i\in\underline{k-1},\  d_i(I)\subset\left[0,\frac{1}{2}\right]\times I^{N-1}; \\\forall j\in\underline n,\ d_j(I)\subset\left[\frac{1}{2},1\right]\times I^{N-1}.\end{matrix*}\right.\right\}.
$$
Then $P^N_k$ and $Q^N_k$ are deformation retracts of $q_k^{-1}([d_{\underline k},x^{\underline k}])$ and $p_o^{N\ -1}(f_k([d_{\underline k},x^{\underline k}]))$ respectively, and the restriction $g_k:P^N_k\rightarrow Q^N_k$ is a fibration with contractible fiber, and therefore a weak equivalence. This implies $g_k:q_k^{-1}([d_{\underline k},x^{\underline k}])\rightarrow p_o^{-1}(f_k([d_{\underline k},x^{\underline k}]))$ is also a weak equivalence.
\end{enumerate}

Therefore by proposition \ref{CondQsiFib2} $p^N_o$ is a quasifibration. That $C^NX_c$ is the fiber of $p_o^N$ follows easily from the definitions.$\blacksquare$

\begin{cor}\label{RelApprox}
    Let $(X_o,X_c)\in\mathcal{SC}^N[\texttt{Top}]$ with $X_o$ well-pointed. Then ${\alpha^N_2}_{(X_o,X_c)}$ is a weak equivalence if $(X_o,X_c)$ are connected or $2<N\leq \infty$ and $(X_o,X_c)$ are grouplike.
\end{cor}

\textbf{Proof:} We have the commutative diagram:
$$
        \xymatrix{
           &C^N X_c\ar[dl]_{\alpha^N}\ar[dd]^{\alpha^N_c}\ar[r]&\pi_o SC^N(X_o,X_c)\ar[r]^{p_o^N}\ar[dd]^{\alpha^N_o}&C^{N-1} X_o\ar[dd]^{\alpha^{N-1}}\\
           \Omega^N\Sigma^NX_c\ar[dr]_\sim&&&\\
           &\Omega^{N}(\Sigma^{N-1}(X_o\wedge I)\vee\Sigma^NX_c)\ar[r]&\Omega^N_\text{rel}\Sigma^N_2(X_o,X_c)\ar@{->>}[r]_{\partial}&\Omega^{N-1}\Sigma^{N-1}X_o
        }
$$
If $2<N\leq\infty$ this is a commutative diagram of $\mathcal H$-maps between admissible $\mathcal H$-spaces. If $(X_o,X_c)$ are grouplike $C^NX_c$ and $C^{N-1}X_o$ are also grouplike.

The bottom horizontal maps form a fibration sequence. By the above theorem the upper horizontal maps form a quasifibration sequence. The corollary then follows from \ref{CorApprox} and the five lemma.$\blacksquare$

\subsection{Compatibility of geometric realization, $\Sigma^N_2$, $\Omega^N_2$ and $SC^N$}

We have the following compatibility results, which are analogous to results in \cite[12]{Ma72}. We omit the proofs since they are easy adaptations of the arguments given there.

\begin{prop}
    The natural maps 
    \begin{align*}
    \tau:|\Sigma^N_{2\bullet}(X_o,X_c)_\bullet|&\rightarrow \Sigma^N_{2}|(X_o,X_c)_\bullet|\\
    [[x,s],u]&\mapsto [[x,u],s].
\end{align*}
are homeomorphisms.
\end{prop}

\begin{prop}\label{G12.2}
    Let $\mathcal O$ be any topological 2-operad and $O$ its associated monad in $\texttt{Top}_\ast^{\{o,c\}}$. Then there is a natural homeomorphism
    \begin{align*}
    \nu:|O_\bullet (X_o,X_c)_\bullet|&\rightarrow O|(X_o,X_c)_\bullet|\\
    [[o_{\underline m},x^{\underline m}],u]&\mapsto\left[o_{\underline m},\prod_{i\in\underline m}[x^i,u]\right]
\end{align*}
    such that the following diagrams are commutative
    $$
    \xymatrix{|(X_o,X_c)_\bullet|\ar[r]^{|\eta_\bullet|}\ar[dr]_\eta&|O_\bullet (X_o,X_c)_\bullet|\ar[d]^\nu\\
    &O|(X_o,X_c)_\bullet|}\ \ \ \ \ \ \ \ \ \ \  \xymatrix{|O^2_\bullet (X_o,X_c)_\bullet|\ar[d]_{|\mu_\bullet|}\ar[r]^{\nu^2}&O^2|(X_o,X_c)_\bullet|\ar[d]^\mu\\
    |O_\bullet (X_o,X_c)_\bullet|\ar[r]_\nu&O|(X_o,X_c)_\bullet|}
    $$
    
    If $((X_o,X_c)_\bullet,\xi_\bullet)\in O[\texttt{Top}]^{\Delta^\text{op}}$, then $(|X_\bullet|,|\xi_\bullet|\nu^{-1})\in O[\texttt{Top}]$ and geometric realization therefore defines a functor $O[\texttt{Top}]^{\Delta^\text{op}}\rightarrow O[\texttt{Top}]$.
\end{prop}

\begin{prop}\label{|loop|}
    For $\iota_{\bullet}:B_{\bullet}\rightarrow Y_{\bullet}\in(\texttt{Top}_\ast^\shortrightarrow)^{\Delta^\text{op}}$ a proper simplicial relative space with each $\iota_{q}$ an $N-1+q$-connected relative space, the map
    \begin{align*}
        \gamma^N_2:|\Omega^N_{2\bullet}\iota_{\bullet}|&\rightarrow \Omega^N_2|\iota_{\bullet}|\\
        [(\alpha,\gamma),u]&\mapsto\left(s\mapsto [(\alpha(s),s'\mapsto\gamma(s)(s')),u]\right)\\
        [\delta,u]&\mapsto\left(t\mapsto [\delta(t),u]\right)
    \end{align*}
    is a weak equivalence.
    
    For $\iota_{\bullet,\bullet}:B_{\bullet,\bullet}\rightarrow Y_{\bullet,\bullet}\in(\texttt{Sp}^\shortrightarrow)^{\Delta^\text{op}}$ a proper simplicial relative spectrum with each $\iota_{p,q}$ a $p+q$-connected relative space, the map $\gamma^\infty_2:|\Omega^\infty_{2\bullet}\iota_{\bullet,\bullet}|\rightarrow \Omega^\infty_2|\iota_{\bullet,\bullet}|$ induced by the $\gamma^N_2$ above is a weak equivalence.
\end{prop}

\begin{prop}\label{|loopAlg|}
    For $\iota_{\bullet}:B_{\bullet}\rightarrow Y_{\bullet}\in(\texttt{Top}_{\ast}^\shortrightarrow)^{\Delta^\text{op}}$ (or $\iota_{\bullet,\bullet}:B_{\bullet,\bullet}\rightarrow Y_{\bullet,\bullet}\in(\texttt{Sp}^\shortrightarrow)^{\Delta^\text{op}}$ if $N=\infty$) the map $\gamma_2^N$ is a $\mathcal{SC}^N$-map. Also, for $(X_o,X_c)_\bullet\in(\texttt{Top}^{2}_\ast)^{\Delta^\text{op}}$ the following diagram is commutative:
    $$
        \xymatrix{
            |SC^N_\bullet (X_o,X_c)_\bullet|\ar[r]^{\nu}\ar[d]_{|\alpha^N_{2\bullet}|}&SC^N| (X_o,X_c)_\bullet|\ar[d]^{\alpha^N_2}\\
            |\Omega^N_{2\bullet}\Sigma^N_{2\bullet}(X_o,X_c)_\bullet|\ar[r]_{\tau\gamma_2^N}&\Omega^N_2\Sigma^N_2|(X_o,X_c)_\bullet|
        }
    $$
\end{prop}

\subsection{Finite relative recognition principle}

The next definition can be thought of as a relative delooping functor.

\begin{defin}
Let $\pi_2^N:\mathcal D^N\rightarrow\mathcal{SC}^N$ be a weak equivalence of $2\mathbb S$-free 2-operads. The classifying space functor for $\mathcal D^N$-spaces is
\begin{align*}
    B_2^N:\mathcal D^N[\texttt{Top}]&\rightarrow\texttt{Top}_\ast^{\shortrightarrow}\\
    (X_o,X_c)&\mapsto B(\Sigma_2^N,D^N,(X_o,X_c))
\end{align*}
\end{defin}

The following is the central theorem of the relative recognition principle.

\begin{teor}\label{FinRelRecogPrin}
    Let $\pi^N_2:\mathcal D^N\rightarrow\mathcal{SC}^N$ be a weak equivalence of $2\mathbb S$-free 2-operads. Let $(X_o,X_c)$ be a $\mathcal D^N$-space such that $(X_o,X_c)$ are well-pointed and consider the following $\mathcal D^N$-maps:
    $$
        \xymatrix@C=2cm{
            B(D^N,D^N,(X_o,X_c))\ar[d]_{\varepsilon(\xi)}\ar[r]^{B(\alpha_2^N\pi^N_2,1,1)} & B(\Omega_2^N\Sigma_2^N,D^N,(X_o,X_c))\ar[d]^{\gamma^N_2}\\
            (X_o,X_c) & \Omega_2^NB_2^N(X_o,X_c)
        }
    $$
    \begin{enumerate}[i)]
        \item $\varepsilon(\xi)$ is a strong deformation retraction with right inverse $\tau(\zeta)$, where $\zeta:(X_o,X_c)\rightarrow D^N(X_o,X_c)$ is given by the unit $\zeta$ of $D^N$;
        
        \item $B(\alpha_2^N\pi^N_2,1,1)$ is a weak equivalence if $X_o$ and $X_c$ are connected or if $2<N<\infty$ and $X_o$ and $X_c$ are grouplike;
        
        \item $\gamma_2^N$ is a weak equivalence;
        
        \item The composite $\gamma_2^N
        B(\alpha_2^N\pi,1,1) \tau(\zeta)$ coincides with $\Omega^N_2(\tau(1_{\Sigma_2^N}))\eta^N_{2}$, and is a weak equivalence if $X_o$ and $X_c$ are connected or if $2<N<\infty$ and $X_o$ and $X_c$ are grouplike;.
        
        \item $B(\Sigma_2^N,D^N,(X_o,X_c))$ is $m+N$-connected if $(X_o,X_c)$ are both $m$-connected.
        
        \item For $(\iota:B\rightarrow Y)\in\texttt{Top}_\ast^\shortrightarrow$, $\varepsilon(\epsilon_2^N):B_2^N\Omega_2^N \iota\rightarrow \iota$ is a weak equivalence if $\iota$ is $N$-connected or if $N>2$ and $\iota$ is $(N-1)$-connected; for all $\iota$, the following diagram is commutative and $\Omega_2^N\varepsilon(\epsilon_2^N)$ is a retraction with right inverse $\Omega^N_2(\tau(1_{\Sigma_2^N\Omega_2^N}))\eta^N_{2\Omega_2^N}$,
        $$
        \xymatrix@C=2cm{
        B(D^N,D^N,\Omega_2^N \iota)\ar[r]^(0.47){B(\alpha_2^N\pi_2^N,1,1)}\ar[d]_{\varepsilon(\theta_2^N\pi_2^N)}&B(\Omega_2^N\Sigma_2^N,D^N,\Omega_2^N \iota)\ar[d]^{\gamma_2^N}\ar[dl]^{\varepsilon(\Omega^N_2\epsilon^N_2)}\\
        \Omega_2^N \iota&\Omega_2^NB_2^N\Omega_2^N \iota\ar[l]^(0.55){\Omega_2^N\varepsilon(\epsilon_2^N)}}
        $$
        
        \item For $(Y_o,Y_c)\in\texttt{Top}_\ast^2$, $\varepsilon\left(\epsilon^N_{2\Sigma^N_2}\Sigma^N_2(\alpha^N_2\pi_2^N)\right):B_2^ND^N(Y_o,Y_c)\rightarrow \Sigma_2^N(Y_o,Y_c)$ is a strong deformation retraction with right inverse $\tau(\Sigma_2^N\zeta)$.
    \end{enumerate}
\end{teor}

\textbf{Proof:} $\varepsilon(\xi)$ and $B(\alpha^N_2\pi_2^N,1,1)$ are realizations of maps of simplicial $\mathcal D^N$-spaces and are therefore maps of $\mathcal D^N$-spaces by \ref{G12.2}. $\gamma^N_2$ is a map of $\mathcal D^N$-spaces by \ref{|loopAlg|}. $i)$ and $vii)$ hold before realization by \cite[9.10, 9.11]{Ma72}, hence after realization by \cite[11.10]{Ma72}. $ii)$ holds before realization by \ref{Act M.3.4} and \ref{RelApprox}, and therefore after realization by \cite[11.13]{Ma72} and \cite[A.1,A.5]{Ma74}. $iii)$ follows from \ref{|loop|}. $(iv)$ follows from $i)$, $ii)$ and $iii)$. $v)$ follows from \cite[11.12]{Ma72} and \cite[A.5]{Ma74}. The upper triangle of $vi)$ commutes by the naturality of $\varepsilon$, and the bottom triangle by \cite[9.11]{Ma72}. The fact that $\varepsilon(\epsilon_2^N)$ is a weak homotopy equivalence when the connectivity assumptions on $\iota$ are met follows from the commutativity of the diagram.$\blacksquare$

The following corollary tells us that for connected (grouplike if $N>2$) $\mathcal D^N$-algebras then its delooping is unique up to weak homotopy equivalence among $N$-connected ($N-1$-connected if $N>2$) relative spaces. 

\begin{cor}
    Under the hypothesis of the theorem, consider the following diagram of $\mathcal D^N$-algebras:
    $$
        \xymatrix{(X_o,X_c)&\ar[l]_{f}(X'_o,X'_c)\ar[r]^(0.6){g}&\Omega^N_2\iota}
    $$
    if $f$ and $g$ are weak equivalences and either $N>2$, all algebras are grouplike and $\iota$ is $N-1$-connected or $N\geq1$, all algebras are connected and $\iota$ is $N$-connected, then the diagram
    $$
        \xymatrix@C=1.5cm{B^N_2(X_o,X_c)&\ar[l]_{B(1,1,f)}B^N_2(X'_o,X'_c)\ar[r]^(0.65
        ){\varepsilon(\epsilon^N_2)B(1,1,g)}&\iota}
    $$
    displays a weak equivalence between $\iota$ and $B^N_2(X_o,X_c)$.
\end{cor}

\textbf{Proof:} $\varepsilon(\epsilon^N_2)$ is a weak equivalence by item $vi)$ of the theorem, and $B(1,1,f)$ and $B(1,1,g)$ are weak equivalences before realization by the approximation theorem, and also after realization by \cite[11.13]{Ma72}.$\blacksquare$

These results can be put together into an equivalence of homotopy categories when $\mathcal D^N$ is cofibrant.

\begin{teor}
    For $2<N<\infty$ and $\pi:\overline{\mathcal{SC}^N}\rightarrow \mathcal{SC}^N$ a cofibrant replacement of the Swiss-cheese operad the functors $B_2^N$ and $\Omega_2^N$ induce equivalences of categories.
    $$(\mathbb LB_2^N\dashv\mathbb R \Omega_2^N):\mathcal Ho \texttt{Top}^{\shortrightarrow}_{N-1}\rightleftharpoons\mathcal Ho \overline{\mathcal{SC}^N}[\texttt{Top}]_\text{grp}.
    $$
    
    For all $N$ they induce an equivalence
    $$(\mathbb LB_2^N\dashv\mathbb R \Omega_2^N):\mathcal Ho \texttt{Top}^{\shortrightarrow}_{N}\rightleftharpoons\mathcal Ho \overline{\mathcal{SC}^N}[\texttt{Top}]_\text{0}.
    $$
\end{teor}

\textbf{Proof:} We prove the first statement, the second follows from similar arguments. Note that by \ref{LeftDerB} the left derived functor $\mathbb LB^N_2$ is well defined, and since all objects in $\texttt{Top}^{\shortrightarrow}$ are fibrant then $\mathbb LB^N_2=B^N_2$. Let $\iota\in\mathcal Ho\texttt{Top}_{N-1}^{\shortrightarrow}$. We define the counit of the adjunction $[\epsilon]$ as the homotopy class of the composition:
    
    $$
        B^N_2\text{Cof}\Omega^N_2\iota\xrightarrow{B(1,1,\text{cof}_{\Omega^N_2\iota})}B^N_2\Omega^N_2\iota\xrightarrow{\varepsilon(\epsilon_2^N)}\iota
    $$
    The first map is a weak equivalence for every $\iota\in\texttt{Top}^\shortrightarrow_{N-1}$ and the second map is a weak equivalence by \ref{FinRelRecogPrin}$.vi)$ because $\iota$ is assumed to be $N-1$-connected.
    
    Let $(X_o,X_c)\in\mathcal Ho \overline{\mathcal{SC}^N}[\texttt{Top}]_{gr}$. Since $B(\overline{SC^N},\overline{SC^N},(X_o,X_c))$ is also cofibrant by the Whitehead theorem there is a homotopy inverse $\varepsilon(\xi)^{-1}$ of $\varepsilon(\xi)$ in $\overline{\mathcal{SC}^N}[\texttt{Top}]$, which is unique up to homotopy and therefore it is homotopy equivalent to $\tau(\zeta)$. Therefore we have the commutative diagram
    $$
        \xymatrix@C=2cm{
            (\ast_o,\ast_c) \ar@{^{(}->}[d]\ar@{^{(}->}[r]& \text{Cof}\Omega^N_2B^N_2(X_o,X_c)\ar@{->>}[d]_\sim^{\text{cof}_{\Omega^N_2B^N_2(X_o,X_c)}}\\
            (X_o,X_c)\ar@{.>}[ur]^{\eta} \ar[r]_(0.45){\gamma^N_2B(\alpha^N_2\pi_2^N,1,1)\varepsilon(\xi)^{-1}} & \Omega^N_2B^N_2(X_o,X_c)
        }
    $$
    with the bottom map a weak equivalence by \ref{FinRelRecogPrin}$.iv)$. We define the unit of the adjunction $[\eta]$ as the homotopy class of the lift of this diagram.

The unit-counit equations hold by the naturality of the cofibrant replacement functor, the homotopy commutativity by \ref{FinRelRecogPrin}$.iv)$ of the upper left and left triangle of the following diagrams respectively, the commutativity of the rest of the diagrams and the $(\Omega^N_2\dashv\Sigma^N_2)$-adjunction.
$$
    \xymatrix@C=3cm{
    B^N_2(X_o,X_c)\ar[dr]^{B^N_2\eta^N_2}\ar[r]^(0.45){B^N_2(\gamma^N_2B(\alpha^N_2\pi^N_2,1,1)\varepsilon(\xi)^{-1})}&B^N_2\Omega^N_2B^N_2(X_o,X_c)\ar[r]^{\varepsilon(\epsilon^N_{2B^N_2})}&B^N_2(X_o,X_c)\\
    &B^N_2\Omega^N_2\Sigma^N_2(X_o,X_c)\ar[rd]^{\varepsilon(\epsilon^N_{2\Sigma^N_2})}\ar[u]_{B^N_2\Omega^N_2(\tau(1_{\Sigma^N_2}))}&\\
    \Sigma^N_2(X_o,X_c)\ar[uu]^{\tau(1_{\Sigma^N_2})} \ar[r]_{\Sigma^N_2\eta^N_2}& \Sigma^N_2\Omega^N_2\Sigma^N_2(X_o,X_c)\ar[u]^{\tau(1_{\Sigma^N_2\Omega^N_2\Sigma^N_2})}\ar[r]_{\epsilon^N_{2\Sigma^N_2}} &\Sigma^N_2(X_o,X_c)\ar[uu]_{\tau(1_{\Sigma^N_2})}
    }
$$
$$
    \xymatrix@C=3cm{
    \Omega^N_2\iota\ar[dr]_{\eta^N_{2\Omega^N_2}}\ar[r]^(0.45){\gamma^N_2B(\alpha^N_2\pi^N_2,1,1)\varepsilon(\theta^N_2\pi^N_2)^{-1}}&\Omega^N_2B^N_2\Omega^N_2\iota\ar[r]^{\Omega^N_2\varepsilon(\epsilon^N_{2})}&\Omega^N_2\iota\\
    &\Omega^N_2\Sigma^N_2\Omega^N_2\iota\ar[ur]_{\Omega^N_2\epsilon^N_{2}}\ar[u]_{\Omega^N_2(\tau(1_{\Sigma^N_2\Omega^N_2}))}&
    }\blacksquare
$$

The statements in this section for $N=1$ can be extended to the non-connected grouplike cases by the relative recognition result in \cite{HLS16}. The difficulty of making the generalization for $N=2$ comes from the fact that $\mathcal C^1$-algebras aren't generally admissible, so $\alpha^2_2$ may not be a weak equivalence on grouplike $\overline{\mathcal{SC}^2}$-spaces, and the proof in \cite{HLS16} uses the fact that $\mathcal{SC}^1$ is weakly equivalent to the 2-opereads of actions of topological monoids on spaces, which is clearly not true for $\mathcal{SC}^2$. Therefore the general case for $N=2$ remains open.

\subsection{Infinite relative recognition principle}

Consider $\pi^N_2:\mathcal{D}^N\rightarrow \mathcal{SC}^N$ for $1\leq N<\infty$ a choice of weak equivalences of $2\mathbb S$-free 2-operads such that there are 2-operad maps $\sigma^N:\mathcal{D}^N\rightarrow\mathcal{D}^{N+1}$ that commute with the inclusions $\sigma^N:\mathcal{SC}^N\rightarrow\mathcal{SC}^{N+1}$. Define $\mathcal D^\infty:=\varinjlim\mathcal D^\infty$. This defines a weak equivalences of $2\mathbb S$-free 2-operads $\pi^\infty_2:\mathcal{D}^\infty\rightarrow \mathcal{SC}^\infty$. Note that any $\mathcal{D}^\infty$-space is a $\mathcal{D}^N$-space for all $1\leq N<\infty$ by the inclusions, and that the homeomorphisms $\mathds S^N\wedge \mathds S^1\rightarrow \mathds S^{N+1}$ and the $\sigma^N$ induce maps 
$$
    \Sigma^{1\shortrightarrow}B(\Sigma_2^{N},\mathcal{D}^{N},(X_o,X_c))\rightarrow B(\Sigma_2^{N+1},\mathcal{D}^{N+1},(X_o,X_c))
$$
which gives this collection of spaces the structure of a relative spectrum. We can then think of the next definition as an infinite relative delooping functor.

\begin{defin}
The classifying space functor for $\mathcal D^\infty$-spaces is defined as
\begin{align*}
    B^\infty_2:\mathcal D^\infty[\texttt{Top}]&\rightarrow\texttt{Sp}^{\shortrightarrow}\\
    (X_o,X_c)&\mapsto B(\Sigma_2^{\bullet+1},\mathcal{D}^{\bullet+1},(X_o,X_c))
\end{align*}
\end{defin}

The following tells us that the relative delooping of grouplike $\mathcal D^\infty$-alagebras are relative $\Omega$-spectra.

\begin{lema}\label{DeloopFib}
    If $(X_o,X_c)\in\mathcal D^\infty[\texttt{Top}]$ is grouplike, then $B^\infty_2(X_o,X_c)$ is fibrant.
\end{lema}

\textbf{Proof:} This follows from the recognition principle in the previous section since each $B^\infty_2(X_o,X_c)_\bullet$ is the $\bullet+1$-fold relative delooping of $(X_o,X_c)$.$\blacksquare$

The results of the previous sections then induce the following.

\begin{teor}\label{InfinRelRecogPrin}
    Let $(X_o,X_c)$ be a $\mathcal D^\infty$-space such that $(X_o,X_c)$ are well-pointed and consider the following morphisms of $\mathcal D^\infty$-spaces:
    $$
        \xymatrix@C=2cm{
            B(D^\infty,D^\infty,(X_o,X_c))\ar[d]_{\varepsilon(\xi)}\ar[r]^{B(\alpha_2^\infty\pi^\infty_2,1,1)} & B(\Omega_2^\infty\Sigma_2^\infty,D^\infty,(X_o,X_c))\ar[d]^{\gamma^N_2}\\
            (X_o,X_c) & \Omega_2^\infty B_2^\infty(X_o,X_c)
        }
    $$
    \begin{enumerate}[i)]
        \item $\varepsilon(\xi)$ is a strong deformation retraction with right inverse $\tau(\zeta)$, where $\zeta:(X_o,X_c)\rightarrow D^\infty(X_o,X_c)$ is given by the unit $\zeta$ of $\mathcal D^\infty$;
        
        \item $B(\alpha_2^\infty\pi^\infty_2,1,1)$ is a weak equivalence if $X_o$ and $X_c$ are grouplike;
        
        \item $\gamma_2^\infty$ is a weak equivalence;
        
        \item The composite $\gamma_2^\infty
        B(\alpha_2^\infty\pi^\infty_2,1,1) \tau(\zeta):(X_o,X_c)\rightarrow \Omega_2^\infty B_2^\infty(X_o,X_c)$ coincides with $\Omega^\infty_2(\tau(1_{\Sigma_2^\infty}))\eta^\infty_{2}$, and is a weak equivalence if $X_o$ and $X_c$ are grouplike;
        
        \item $B(\Sigma_2^\infty,D^\infty,(X_o,X_c))$ is $m+1$-connective if $(X_o,X_c)$ are both $m$-connected;
        
        \item For $(\iota_\bullet:B_{\bullet}\rightarrow Y_{\bullet})\in\texttt{Sp}^\shortrightarrow$, $\varepsilon(\epsilon_2^\infty):B_2^\infty\Omega_2^\infty \iota_\bullet\rightarrow \iota_\bullet$ is a weak equivalence if $\iota_\bullet$ is $0$-connective; for all $\iota_\bullet$, the following diagram is commutative and $\Omega_2^\infty\varepsilon(\epsilon_2^\infty)$ is a retraction with right inverse $\Omega^\infty_2(\tau(1_{\Sigma_2^\infty\Omega_2^\infty}))\eta^\infty_{2\Omega_2^\infty}$,
        $$
        \xymatrix@C=2cm{
        B(D^\infty,D^\infty,\Omega_2^\infty \iota_\bullet)\ar[r]^(0.47){B(\alpha_2^\infty\pi_2^\infty,1,1)}\ar[d]_{\varepsilon(\theta_2^\infty\pi^\infty_2)}&B(\Omega_2^\infty\Sigma_2^\infty,D^\infty,\Omega_2^\infty \iota_\bullet)\ar[d]^{\gamma_2^\infty}\ar[dl]^{\varepsilon(\Omega^\infty_2\epsilon^\infty_2)}\\
        \Omega_2^\infty \iota_\bullet&\Omega_2^\infty B_2^\infty\Omega_2^\infty \iota_\bullet\ar[l]^(0.55){\Omega_2^\infty\varepsilon(\epsilon_2^\infty)}}
        $$
        
        \item For $(Y_o,Y_c)\in\texttt{Top}_\ast^2$, $\varepsilon\left(\epsilon^\infty_{2\Sigma^\infty_2}\Sigma^\infty_2(\alpha^\infty_2\pi_2^\infty)\right):B_2^\infty D^\infty(Y_o,Y_c)\rightarrow \Sigma_2^\infty(Y_o,Y_c)$ is a strong deformation retraction with right inverse $\tau(\Sigma_2^\infty\zeta)$.
    \end{enumerate}
\end{teor}

\begin{cor}
    Under the hypothesis of the theorem, consider the following diagram:
    $$
        \xymatrix{(X_o,X_c)&\ar[l]_{f}(X'_o,X'_c)\ar[r]^(0.6){g}&\Omega^\infty_2\iota_\bullet}
    $$
    if $f$ and $g$ are weak equivalences of group-like $\mathcal D^\infty$-algebras and $\iota_\bullet$ is $0$-connective, then the diagram
    $$
        \xymatrix@C=1.5cm{B^\infty_2(X_o,X_c)&\ar[l]_{B(1,1,f)}B^\infty_2(X'_o,X'_c)\ar[r]^(0.65
        ){\varepsilon(\epsilon^N_2)B(1,1,g)}&\iota_\bullet}
    $$
    displays a weak equivalence between $\iota$ and $B^N_2(X_o,X_c)$.
\end{cor}

In the following equivalence of categeories, by lemma \ref{DeloopFib}, we can again use $\mathbb L B^\infty_2=B^\infty_2$.

\begin{teor}
    If $\pi_2^N:\overline{\mathcal{SC}^N}\rightarrow \mathcal{SC}^N$ are cofibrant replacements satisfying the conditions on the begining of this section, the functors $B_2^\infty$ and $\Omega_2^\infty$ induce an equivalence of categories.
    $$(\mathbb LB_2^\infty\dashv\mathbb R \Omega_2^\infty):\mathcal Ho \texttt{Sp}^{\shortrightarrow}_{0}\rightleftharpoons\mathcal Ho \overline{\mathcal{SC}^\infty}[\texttt{Top}]_\text{grp}.
    $$
\end{teor}


\begin{thebibliography}{56}

\bibitem[BF78]{BF78} Bousfield, Aldridge; Friedlander, Eric. \emph{Homotopy theory of $\Gamma$-spaces, spectra, and bisimplicial sets.} Springer Lecture Notes in Math., 658 (1978): 80-130.

\bibitem[BM03]{BM03} Berger, Clemens; Moerdijk, Ieke. \emph{Axiomatic homotopy theory for operads.} Commentarii Mathematici Helvetici 78.4 (2003): 805-831.

\bibitem[BM07]{BM07} Berger, Clemens; Moerdijk, Ieke. \emph{Resolution of coloured operads and rectification of homotopy algebras.} Contemporary Mathematics 431 (2007): 31-58.

\bibitem[BV68]{BV68} Boardman, John Michael; Vogt, Rainer M. \emph{Homotopy-everything H-spaces.} Bulletin of the American mathematical society 74.6 (1968): 1117-1122.

\bibitem[CLM76]{CLM76} Cohen, F. Ronald; Lada, T. Joseph; May, J. Peter. \emph{The homology of iterated loop spaces}. Lecture Notes in Mathematics, 533 (1976).

\bibitem[ESt67]{ESt67} Steenrod, N. E. \emph{A convenient category of topological spaces.} Michigan Math. J., 14 no. 2 (1967): 133--152.

\bibitem[Fr13]{Fr13} Frankhuizen, Robin. \emph{The Recognition Principle for Grouplike $\mathcal{C}_n$-algebras}. MS thesis. 2013.

\bibitem[Hi09]{Hi09} Hirschhorn, Philip S. \emph{Model categories and their localizations.} American Mathematical Soc., 99 (2009).

\bibitem[Hi15]{Hi15} Hirschhorn, Philip S. "The Quillen model category of topological spaces." arXiv preprint arXiv:1508.01942 (2015).

\bibitem[HLS16]{HLS16} Hoefel, Eduardo; Livernet, Muriel; Stasheff, Jim. \emph{$A_\infty$-actions and recognition of relative loop spaces}. Topology and its Applications
206 (2016): 126-147.

\bibitem[HL12]{HL12} Hoefel, Eduardo; Livernet, Muriel. \emph{Open-closed homotopy algebras and strong homotopy Leibniz pairs through Koszul operad theory.} Letters in Mathematical Physics 101.2 (2012): 195-222.

\bibitem[HoE07]{HoE07} Hoefel, Eduardo. \emph{OCHA and the Swiss-cheese operad.} arXiv preprint arXiv:0710.3546 (2007).

\bibitem[HoM07]{HoM07} Hovey, Mark. \emph{Model categories.} American Mathematical Soc., 63 (2007).

\bibitem[Id17]{Id17} Idrissi, Najib. \emph{Swiss-cheese operad and Drinfeld center.} Israel Journal of Mathematics 221.2 (2017): 941-972.

\bibitem[Ko99]{Ko99} Kontsevich, Maxim. \emph{Operads and motives in deformation quantization.} Letters in Mathematical Physics 48.1 (1999): 35-72.

\bibitem[KS06]{KS06} Kajiura, Hiroshige; Stasheff, Jim. \emph{Homotopy algebras inspired by classical open-closed string field theory.} Communications in mathematical physics 263.3 (2006): 553-581.

\bibitem[Li15]{Li15} Livernet, Muriel. \emph{Non-formality of the Swiss-cheese operad.} Journal of Topology 8.4 (2015): 1156-1166.

\bibitem[Lu09]{Lu09} Lurie, Jacob. \emph{Higher Topos Theory.} Princeton University Press AM-170 (2009).

\bibitem[Ma72]{Ma72} May, J. Peter. \emph{The geometry of iterated loop spaces}, Lecture Notes in Mathematics, 271 (1972).

\bibitem[Ma74]{Ma74} May, J. Peter. \emph{$E_\infty$-spaces, group completions and permutative categories.} New Developments in Topology, London Math. Soc. Lecture Note Series, 11 (1974).

\bibitem[Ma90]{Ma90} May, J. Peter. \emph{Weak equivalences and quasifibrations.} Lecture Notes in Mathematics, 1425 (1990): 91-101.

\bibitem[Ma99]{Ma99} May, J. Peter. \emph{A concise course in algebraic topology.} University of Chicago press, (1999).

\bibitem[Ma09]{Ma09} May, J. Peter. \emph{What precisely are $E_{\infty}$ ring spaces and $E_{\infty}$ ring spectra.} New topological contexts for Galois theory and algebraic geometry (BIRS 2008) 16 (2009): 215-282.

\bibitem[MDS76]{MDS76} McDuff, Dusa; Segal, Graeme. \emph{Homology fibrations and the “group-completion” theorem.} Inventiones mathematicae, 31.3 (1976): 279-284.

\bibitem[Qu15]{Qu15} Quesney, Alexandre. \emph{Swiss Cheese type operads and models for relative loop spaces.} arXiv preprint arXiv:1511.05826 (2015).

\bibitem[Ri09]{Ri09} Riehl, Emily. "A concise definition of a model category." Preprint, available at http://www.math.harvard.edu/$\sim$eriehl/modelcat.pdf (2009).

\bibitem[StJ63]{StJ63} Stasheff, James D. \emph{Homotopy associativity of H-spaces I.} Trans. Amer. Math. Soc. 108 (1963), 275-292.

\bibitem[StN09]{StN09} Strickland, Neil P. \emph{The category of CGWH spaces.} preprint, 12 (2009).

\bibitem[Sc97]{Sc97} Schwede, Stefan. \emph{Spectra in model categories and applications
to the algebraic cotangent complex}. Journal of Pure and Applied Algebra, 120 (1997): 77-104.

\bibitem[Vo99]{Vo99} Voronov, Alexander A. \emph{The Swiss-cheese operad.} Contemporary Mathematics 239 (1999): 365-374.

\bibitem[Wi17]{Wi17} Willwacher, Thomas. \emph{(Non-) formality of the extended Swiss Cheese operads.} arXiv preprint arXiv:1706.02945 (2017).

\bibitem[Zw98]{Zw98} Zwiebach, Barton. \emph{Oriented open-closed string theory revisited.} Annals of Physics 267.2 (1998): 193-248.
\end{thebibliography}
\end{document}